\documentclass[a4paper,oneside,reqno,svgnames]{amsart}
\usepackage[utf8x]{inputenc}
\usepackage[english]{babel}
\usepackage{lmodern} 

\usepackage{geometry}

\usepackage{amsmath,amsfonts,amssymb,amsthm}
\usepackage{mathtools}
\usepackage{bm} 			
\usepackage{stmaryrd} 		
\usepackage{bbm} 			
\usepackage{mathrsfs}   		
\usepackage{dsfont}  		

\usepackage{float}
\usepackage{tikz-cd}

\usepackage{enumitem} 		

\usepackage{xcolor}
\usepackage[colorlinks, linkcolor=DarkBlue, citecolor=cyan, backref=page]{hyperref} 
\usepackage{csquotes} 
\usepackage[textsize=tiny,textwidth=1.08in]{todonotes}
\usepackage[bottom]{footmisc}

\theoremstyle{plain}
\newtheorem{theorem}{Theorem}[section]
\newtheorem{proposition}[theorem]{Proposition}
\newtheorem{lemma}[theorem]{Lemma}
\newtheorem{corollary}[theorem]{Corollary}

\newtheorem*{theorem*}{Theorem}
\newtheorem*{proposition*}{Proposition}
\newtheorem*{lemma*}{Lemma}
\newtheorem*{corollary*}{Corollary}
\newtheorem*{property*}{Properties}
\newtheorem*{conjecture*}{Conjecture}

\theoremstyle{definition}
\newtheorem{definition}[theorem]{Definition}

\newtheorem{remark}[theorem]{Remark}
\newtheorem*{definition*}{Definition}
\newtheorem*{example*}{Example}
\newtheorem*{remark*}{Remark}

\theoremstyle{remark}

\newtheorem*{note*}{Remark}
\newtheorem*{exercise*}{Exercise}
\newtheorem*{notation*}{Notation}

\newcommand{\m}[1]{\mathcal{#1}}
\newcommand{\bb}[1]{\mathbb{#1}}
\newcommand{\mbf}[1]{\mathbf{#1}}
\newcommand{\mrm}[1]{\mathrm{#1}}
\newcommand{\f}[1]{\mathfrak{#1}} 
\newcommand{\scr}[1]{\mathscr{#1}}

\newcommand{\idmatrix}{\mathbbm{1}}   
\newcommand{\del}{\partial}         
\newcommand{\delbar}{\bar{\partial}} 
\newcommand{\dd}{\mrm{d}}
\newcommand{\acts}{\curvearrowright} 

\newcommand{\Aut}{\mathrm{Aut}}

\newcommand{\ext}[1]{\bigwedge\nolimits^{#1}}

\newcommand{\Ok}[1]{O\left(k^{#1}\right)}

\DeclarePairedDelimiter{\set}{\{}{\}}   

\DeclareRobustCommand{\SkipTocEntry}[5]{}

\numberwithin{equation}{section} 

\author{Annamaria Ortu}
\title{The analytic moduli space of holomorphic submersions}
\address{SISSA, Via Bonomea 265, 34136 Trieste, Italy}
\email{aortu[at]sissa[dot]it}

\begin{document}

\begin{abstract}
We construct a moduli space that parametrises stable proper holomorphic submersions over a fixed compact K\"ahler base.
Stability is described in terms of the existence of a canonical relatively K\"ahler metric on the submersion, called an optimal symplectic connection. The construction of the moduli space combines techniques from geometric invariant theory with the study of the geometric PDE defining an optimal symplectic connection.
A special case of this moduli space is the moduli space of vector bundles over a compact K\"ahler manifold. We also show that the moduli space is a Hausdorff complex space equipped with a Weil-Petersson type K\"ahler metric.
\end{abstract}

\maketitle
\tableofcontents
\thispagestyle{empty}

\section{Introduction}
A fundamental result in the study of holomorphic vector bundles is the construction of the moduli space of stable vector bundles. Seshadri \cite{Seshadri_ModuliVectorBundlesCurves} and Newstead \cite{Newstead_CharacteristicClassesStableVectorBundles, Newstead_ModuliproblemsOrbitspaces}  gave the first construction of such a moduli space for vector bundles over a curve, while Mumford \cite{Mumford_Stability} introduced semistability in the sense of geometric invariant theory (GIT) to study this moduli space and established its structure as a global GIT quotient.
These moduli spaces also have an analytic interpretation: the Hitchin-Kobayashi correspondence of Donaldson, Uhlenbeck and Yau \cite{Donaldson_HYMstability, UhlenbeckYau_HYMstability} states that a holomorphic vector bundle on a compact K\"ahler manifold is slope-stable if and only if it admits a Hermite-Einstein connection, thus establishing an equivalence between an algebro-geometric stability condition and the solvability of a geometric PDE.
Fujiki and Schumacher \cite{FujikiSchumacher_ModuliHE} later directly constructed the moduli space of Hermite-Einstein vector bundles over a fixed compact K\"ahler manifold using analytic techniques.
These moduli spaces remain a central object of study to this day, and we refer to Greb-Sibley-Toma-Wentworth \cite{GSTW_modulispace_vbd} for recent work containing a discussion of analytic and algebraic compactifications.

Motivated by the Hitchin-Kobayashi correspondence and the construction of the moduli space of stable vector bundles, we extend these results to the context of proper holomorphic submersions.
We recall that a holomorphic vector bundle $F\to B$ induces a proper holomorphic submersion by considering its projectivisation $\bb{P}(F) \to B$.

More precisely, let $\pi_Y:(Y, H_Y) \to (B,L)$ be a proper holomorphic submersion, where $B$ is a smooth projective manifold with an ample line bundle $L$ and $H_Y$ is a relatively ample line bundle over $Y$.
The base $(B,L)$ is considered fixed throughout.
The stability condition we impose on such a submersion is given by a stability condition on the fibres, defined in terms of K-stability, and a global condition on the fibration, given in terms of a canonical choice of a relatively K\"ahler metric, called an \emph{optimal symplectic connection}.
In the following, we consider only polarised K\"ahler classes, although this is not essential. Indeed all our results hold if $c_1(H_X)$ and $c_1(L)$ are replaced respectively by a relative K\"ahler class and a K\"ahler class that do not come from holomorphic line bundles.

When all the fibres admit a constant scalar curvature K\"ahler metric (hence are K-polystable), Dervan and Sektnan \cite{DervanSektnan_OSC1} introduced optimal symplectic connections as solutions to a geometric PDE on the fibration, which produces a canonical relatively K\"ahler metric.
In particular, when the fibration is obtained as the projectivisation of a stable holomorphic vector bundle, an optimal symplectic connection is precisely induced by and induces the Hermite-Einstein connection on the vector bundle \cite[\S3.5]{DervanSektnan_OSC1}.
Similarly, on isotrivial fibrations, the optimal symplectic connection condition becomes the Hermite-Yang-Mills condition on an associated principal bundle \cite{McCarthy_OSC}.

In \cite{Ortu_OSCdeformations} the author generalised the notion of optimal symplectic connections to smooth fibrations whose fibres may not have a constant scalar curvature K\"ahler metric (cscK), but instead are \emph{analytically K-semistable}.
Allowing K-semistable fibres is essential for the construction of the moduli space of fibrations with an optimal symplectic connection. Indeed, when deforming a fibration with cscK fibres, one cannot expect that the fibres remain cscK. Analytic K-semistability, on the other hand, is an open condition, and this allows us to study the local behaviour of families of fibrations with an optimal symplectic connection.
Our main result is the following.

\begin{theorem}
There exists a moduli space $\m{M}$ that parametrises holomorphic submersions over a fixed base, with discrete relative automorphism group and which admit an optimal symplectic connection. The space $\m{M}$ is a Hausdorff complex space and it carries a Weil-Petersson type K\"ahler metric.
\end{theorem}

We now explain briefly the optimal symplectic connection condition.
The assumption that the fibres are analytically K-semistable means that they each admit a degeneration to a cscK manifold. We further assume that these degenerations vary holomorphically in $B$, so that we have a family $(\m{X}, \m{H}) \to B \times S$, parametrised by $S \subseteq \bb{C}$, such that $(\m{X}_s, \m{H}_s) \to B$ is isomorphic to $(Y, H_Y) \to B$ for $s\ne 0$ and at $s=0$ we have a fibration $(X, H_X) \to B$ with cscK fibres.
Applying a relative version of Ehresmann's fibration theorem, $Y$ and $X$ are isomorphic as smooth manifolds. The Chern classes $c_1(H_Y)$ and $c_1(H_X)$, viewed as cohomology classes on the underlying smooth manifold, are equal because they are integral classes. Thus we can take the perspective of fixing the underlying symplectic fibration and deforming the complex structure.
We also fix the complex structure of the base and assume that the deformations keep the projection to $B$ holomorphic.

A relatively symplectic form $\omega$ on $(Y, H_Y)\to (B,L)$ induces a splitting of the tangent bundle of $Y$ into a vertical part, defined as the tangent bundle to the fibres, and a horizontal part, defined using orthogonality with respect to $\omega$. In the language of symplectic geometry, $\omega$ is called a symplectic connection.
A relatively K\"ahler form $\omega$ is an \emph{optimal symplectic connection} on $(Y, H_Y)$ if it satisfies the geometric PDE
\begin{equation}\label{Eq:OSCintro}
p_E(\Delta_{\m{V}}(\Lambda_{\omega_B} (m^*F_{\m{H}})) + \Lambda_{\omega_B} \rho_{\m{H}}) + \frac{\lambda}{2}\nu = 0.
\end{equation}
In this expression $F_{\m{H}}$ and $\rho_{\m{H}}$ are curvature quantities which depend on $\omega$, $\nu$ is a curvature quantity that depends on the infinitesimal change in the complex structure and $\lambda >0$ is a constant. The left-hand side is a smooth function on $Y$, and the map $p_E$ is the projection onto the global sections of the vector bundle $E \to B$ of fibrewise holomorphy potentials with respect to the relatively cscK complex structure of $X$.
In particular, the equation should be interpreted as a second-order elliptic PDE on the vector bundle $E$.

We construct the moduli space $\m{M}$ by gluing local charts around fibrations that admit an optimal symplectic connection. If $(Y,H_Y) \to (B,L)$ is such a fibration, the local moduli space around $Y$ is given as the quotient
\[
W_Y/\Aut(\pi_Y),
\]
where $W_Y$ is the complex space of deformations of $Y$ which also admit an optimal symplectic connection, and we quotient by the action of the discrete group $\Aut(\pi_Y)$ of relative automorphisms, which is finite.

The definition of $W_Y$ essentially involves two steps.
The first step, established in \S \ref{Subsec:openness_setting}, consists of finding a locally closed analytic space which parametrises all small deformations of the complex structure of $Y$ that admit a degeneration to a fibration with cscK fibres.
To do so, we combine the theory of deformations of cscK manifolds of Sz{\'e}kelyhidi \cite{Szekelyhidi_deformations} and Br\"onnle \cite{Bronnle_PhDthesis} and the theory of deformations of fibrations, developed by the author in \cite{Ortu_OSCdeformations}. More precisely, we consider a fibration $(Y, H_Y) \to (B,L)$ with an optimal symplectic connection and a degeneration to a fibration $(X,H_X)\to (B,L)$ with cscK fibres. A fibration version of Kuranishi's theorem allows us to parametrise the compatible vertical deformations of the complex structure of $X$ with a complex space $V_\pi$ of harmonic $(0,1)$-forms with values in the $(1,0)$-vertical tangent bundle.
Working locally in $B$, we then establish that the deformations of $Y$ which degenerate to a relatively cscK fibration form a locally closed analytic subvariety $V_\pi^+$ of $V_\pi$. We explicitly construct the relatively cscK degeneration using techniques from GIT; although we do not use it directly, our construction is related to the Bya{\l}inicki-Birula decomposition \cite{BB_decomposition, BB_decomposition_reductive_groups}. This is the key new step in our construction, not present in other constructions of moduli spaces.

The second step, established in \S \ref{Subsec:openness_OSC}, consists of proving that, in a small open neighbourhood $W_Y$ of the point associated to $Y$ in $V_\pi^+$, all fibrations admit an optimal symplectic connection. The proof of this essentially relies on the implicit function theorem and employs the linearisation of the equation studied in \cite{Ortu_OSCdeformations} with respect to the complex structure, and it is here that we use the assumption on the discreteness of the relative automorphism group.

In the definition of an optimal symplectic connection and in the construction of the space $W_Y$ it is essential to assume that the connected component of the identity of the groups of automorphisms $\Aut_0(X_b, H_b)$ of the fibres of the relatively cscK degeneration are all isomorphic. This assumption is considered to be a fixed datum in our construction of the moduli space.

While our work concentrates on the analytic aspects of the moduli space of holomorphic submersions, it is natural to ask if it is possible to give an algebro-geometric construction based on stability.
In fact, there exists a notion of stability of fibrations developed by Dervan, Sektnan \cite{DervanSektnan_OSC2}, and Hallam \cite{Hallam_geodesics}, which comes from K-stability and is equivalent to slope stability when the fibration is the projectivisation of a vector bundle.
An algebro-geometric construction would lead to the structure of a variety on the moduli space rather than the structure of a complex space, and would naturally allow singular fibres.
Such a moduli space would parametrise certain stable fibrations which degenerate to a fibration whose general fibre is K-polystable. Moreover, the automorphism group of the general fibre of the degeneration should be fixed.
When singular fibres are present, we expect that the stability of the fibration should be interpreted in terms of existence of solutions to the PDE \eqref{Eq:OSCintro} on the smooth locus having a prescribed behaviour near the singularities.

Recently, Hashizume and Hattori \cite{HashizumeHattori_Moduli_CYfibrations} constructed a quasi-projective moduli space which parametrises adiabatically K-stable fibrations over a curve where the generic fibre is Calabi-Yau, using algebro-geometric techniques. Adiabatic K-stability \cite{Hattori_adiabatic_Kstability} is a condition defined on a fibration with K-polystable fibres over a twisted K-polystable base, and is related to Dervan-Sektnan's stability of fibrations and the existence of optimal symplectic connections.
In the special case where the fibres are all smooth, our proof constitutes an alternative construction of Hashizume-Hattori's result.
However, our construction of the moduli space of optimal symplectic connections does not require any condition on the base of the submersion and allows the fibres to have automorphisms.
For example, the projectivisation of a stable simple vector bundle over a K-unstable base is not adiabatically K-stable, yet it admits an optimal symplectic connection, thus it is represented by a point in our moduli space. This observation is based on a discussion with M. Hattori, whom we gratefully thank.

\addtocontents{toc}{\SkipTocEntry}
\subsection*{Outline}
In Section \ref{Sec:preliminaries} we give preliminary definitions and results on deformations of cscK manifolds and we describe the deformation theory of submersions with cscK fibres. We then explain briefly the theory of optimal symplectic connections.
In Section \ref{Sec:analytic_moduli_space} we construct the moduli space of holomorphic submersions admitting an optimal symplectic connection.
Then in Section \ref{Sec:WPmetric} we describe a Weil-Petersson type K\"ahler metric on the moduli space, along with a natural line bundle.

\addtocontents{toc}{\SkipTocEntry}
\subsection*{Acknowledgments}
I am grateful to my supervisors, Ruadha\'i Dervan and Jacopo Stoppa, for many enlightening discussions and for their encouragement. I also thank Eloise Hamilton, John Benjamin McCarthy and G\'abor Sz{\'e}kelyhidi for explaining to me much of the GIT theory contained in Section \ref{Subsec:openness_setting}, and Carlo Scarpa and Lars Martin Sektnan for helpful comments.
Part of this work was conducted while visiting Ruadha\'i Dervan at the University of Cambridge and at the University of Glasgow, supported by his Royal Society University Research Fellowship.

\section{Preliminaries}\label{Sec:preliminaries}
This section recalls some basic definitions and properties of K\"ahler manifolds and introduces the theory of optimal symplectic connections.
Let $(M,L)$ be a compact K\"ahler manifold of complex dimension $n$ with an ample line bundle $L$, and let $\omega \in c_1(L)$ be a K\"ahler form. The \emph{scalar curvature} of $\omega$ is the function on $M$ defined as the contraction of the Ricci curvature:
\begin{equation*}
\mrm{Scal}(\omega) = \Lambda_{\omega} \mrm{Ric}(\omega),
\end{equation*}
where
\begin{equation*}
\mrm{Ric}(\omega) = -\frac{i}{2\pi} \del \delbar \mrm{log} \ \mrm{det} (\omega).
\end{equation*}
We will be interested in K\"ahler metrics with constant scalar curvature, where the constant is a topological constant given by the intersection product
\[
\widehat{S} = \frac{n \ c_1(M) \cdot c_1(L)^{n-1}}{c_1(L)^n}.
\]
Let $g_J$ be the Riemannian metric on $M$ induced by $J$ and $\omega$.
A smooth function $h$ on $M$ is called a \emph{holomorphy potential} if the $(1,0)$-part of the Riemannian gradient of $h$, denoted $\nabla_gh$, is a holomorphic vector field.
Consider the operator $\m{D} : C^\infty(M, \bb{C}) \to \Omega^{0,1}(T^{1,0}M)$ defined by
\[
\m{D}(\varphi) = \delbar \nabla^{1,0}_g \varphi.
\]
The \emph{Lichnerowicz operator} is the operator $\m{D}^*\m{D}$, where the adjoint is defined with respect to the $L^2(g)$-inner product.
The Lichnerowicz operator is a 4th-order elliptic operator, and its kernel, which by ellipticity is the kernel of $\m{D}$, coincides with the space of holomorphy potentials on $M$.

\subsection{Deformation theory of cscK manifolds}\label{Subsec:deformations_cscK}
In this section we describe the deformations of the complex structure of a K\"ahler manifold with constant scalar curvature, following Sz{\'e}kelyhidi \cite[\S3]{Szekelyhidi_deformations}; similar results were obtained also by Br\"onnle \cite[Part 1]{Bronnle_PhDthesis}.

Let $(M, \omega)$ be a symplectic manifold. Define the infinite-dimensional manifold
\[
\scr{J} = \set*{J : TM \to TM \ \text{almost complex structure compatible with }\omega}.
\]
We denote by $\scr{J}^{\mrm{int}}$ the complex subspace of \emph{integrable} almost complex structures. The tangent space at a point $J \in \scr{J}$ can be identified with
\[
T_J^{0,1}\scr{J} = \set*{\left. \alpha \in \Omega^{0,1}(T^{1,0}M)\right\vert \omega(\alpha(u), v) + \omega(u, \alpha(x)) = 0}.
\]

Recall that a vector field $\xi_h$ is Hamiltonian with respect to $\omega$ if there exists $h \in C^\infty(M, \bb{R})$ such that $\omega (\xi_h, \cdot)= -\dd h $.
On a K\"ahler manifold we can characterise a Hamiltonian vector field with Hamiltonian function $h$ as $\xi_h =J\nabla^g(h)$, where $\nabla^gh$ is the Riemannian gradient of $h$.

Let $\scr{G}$ be the group of Hamiltonian symplectomorphisms of $(M, \omega)$. Then $\scr{G}$ acts on $\scr{J}$ by pull-back.
The Lie algebra of $\scr{G}$ can be identified with the space of smooth functions on $M$ with mean value zero, denoted $C_0^\infty(M)$.
The following theorem is due to Fujiki \cite{Fujiki_momentmap} and Donaldson \cite{Donaldson_momentmap}.
 \begin{theorem}
 The action $\scr{G} \acts \scr{J}$ is Hamiltonian with moment map
 \begin{equation}\label{Eq:moment_map_cscK}
 \begin{aligned}
 \m{S}: \scr{J} &\longrightarrow \mrm{Lie}(\scr{G})^* \\
 J &\longmapsto \mrm{Scal}(\omega,J) - \widehat{S}.
 \end{aligned}
 \end{equation}
 \end{theorem}
Therefore cscK metrics on $M$ correspond to $J \in \scr{J}^{\mrm{int}}$ such that $\mu(J) = 0$.
The function $\mrm{Scal}(\omega,J) - \widehat{S}$ is the Hermitian scalar curvature if $J$ is non-integrable, which coincides with the Riemannian scalar curvature in the integrable case.
We view it as an element in $C_0^\infty(M)^*$ by identifying $\mrm{Lie}(\scr{G})$ with its dual.

Let $J$ be an integrable compatible complex structure on $(M,\omega)$, and denote by $\m{K}_J(\omega)$ the set of K\"ahler metrics in the same K\"ahler class of $\omega$ with respect to $J$, i.e\
\[
\m{K}_J(\omega) = \set*{\omega' \in [\omega] \mid \omega' = \omega + i \del_J \delbar_J \phi \ \text{for some} \ \phi \in C^\infty(M, \bb{R})}.
\]
For fixed $J \in \scr{J}$, consider the infinitesimal action of $\scr{G}$ on $\scr{J}$, defined by the operator
\begin{equation}\label{Eq:map_P}
\begin{aligned}
P : C_0^\infty(M, \bb{R}) &\to T_J\scr{J} \\
h &\mapsto \m{L}_{\xi_h} J,
\end{aligned}
\end{equation}
where $\xi_h$ is the Hamiltonian vector field with function $h$; let $P^{\bb{C}}$ be its extension to $C_0^\infty(M, \bb{C})$.
The following proposition (see also \cite[\S6.1]{Szekelyhidi_book}) justifies the fact that instead of moving the K\"ahler metric inside the K\"ahler class with respect to a fixed $J$, one can think of $\omega$ as fixed and move the complex structure.

\begin{proposition}\label{Prop:parallel_cpx-symp_structure}
\cite[p.17]{Donaldson_Symmetric_spaces} For every $\omega_\phi \in \m{K}_J(\omega)$ there exists an $f \in \mrm{Diff}_0(M)$ such that $f^*\omega_\phi = \omega$ and $(M, \omega_\phi, J)$ is isomorphic to $(M, \omega, f^*J)$.
In particular, there is a map $F$ that sends a K\"ahler potential $\phi$ to $F_{\phi}(J) = f^*J$. The differential at $0$ of $F$ is given by
\begin{equation*}
\dd_0F(\phi) = J\m{L}_{\xi_{\phi}(\omega)}J.
\end{equation*}
\end{proposition}
Although a complexification of $\scr{G}$ does not genuinely exist, we can think of the image of the map $F$ as an infinitesimal complexified orbit of the action of $\scr{G}$ on $\scr{J}$. The proposition says that a variation of the K\"ahler form in a given K\"ahler class for $J$ fixed corresponds to a variation of the complex structure $J$ in the same $\scr{G}^c$-orbit, for $\omega$ fixed.

The following result characterises $P$ as the adjoint operator of the differential of the scalar curvature when varying the complex structure.
\begin{lemma}[{\cite[\S4]{Donaldson_momentmap}}]\label{Lemma:PofQ}
Let $Q$ be the differential of the scalar curvature map \eqref{Eq:moment_map_cscK} at a point $J \in \scr{J}$. Then if $\alpha \in T_{J}\scr{J}$ and $\varphi \in C^\infty(M)$,
\[
\langle Q(\alpha), \varphi \rangle_{L^2} = \frac{1}{2}\langle \alpha, P(\varphi) \rangle_J
\]
so $P^* = 2Q$, and
\[
Q(P(\phi)) = \mrm{Re}(\m{D}^*\m{D}\phi).
\]
\end{lemma}

We next let $(M, L)$ be a compact K\"ahler manifold with an ample line bundle $L$, and assume that $\omega \in c_1(L)$ has constant scalar curvature with respect to the complex structure $J$. The deformations of $J$ compatible with $\omega$ can be considered as elements of the vector space 
\begin{equation}\label{Eq:Htilde_1_mfd}
\widetilde{H}^1 = \ker \lbrace\square := PP^*+(\delbar^*\delbar)^2\rbrace
\end{equation}
on $T_J\scr{J}$.
Let $K$ be the stabilizer of $J$ with respect to the action of $\scr{G}$ on $\scr{J}$. Then $K$ is the group of Hamiltonian isometries of $(\omega, J)$.
The complexification $K^{\bb{C}}$ is the connected component of the identity of the group of automorphisms of $M$ which lift to $L$, denoted $\Aut_0(M, L)$. The group $K$ (as well as its complexification) act on $\widetilde{H}^1$ by pull-back, and the origin is a fixed point of the action.

The following theorem is a symplectic version of Kuranishi's theorem \cite{Kuranishi_family_cpx_str} and an infinite-dimensional version of Luna's slice theorem \cite{Luna_slice}: it provides a parametrisation of compatible infinitesimal deformations of $J$. A proof that takes into account the compatibility with the symplectic structure is explained in \cite[Lemma 6.1]{ChenSun_CalabiFlow}.

\begin{theorem}\label{Thm:Kuranishi}
There exists a ball around the origin $V \subset \widetilde{H}^1$ and a $K$-equivariant embedding
\begin{equation}\label{Eq:Kuranishi_map}
\Psi : V \to \scr{J}
\end{equation}
such that $\Psi (0) = J$ and
\begin{enumerate}
\item the $\scr{G}^c$-orbit of every integrable complex structure near $J_0$ intersects the image of $\Psi$;
\item If $x, x' \in V$ are in the same orbit for the complexified action of $K$, and $\Psi(x)$ is integrable, then $\Psi(x), \Psi(x')$ are in the same $\scr{G}^c$-orbit;
\item $\mrm{Scal}(\omega, \Psi(x))$ is an element of the Lie algebra of $K$.
\end{enumerate}
\end{theorem}
In particular, a moment map for the action of $K$ on $V$ is given by
\begin{equation}\label{Eq:moment_map_fd}
x \mapsto \mrm{Scal}(\omega, \Psi(x)) -\widehat{S}.
\end{equation}
The Kuranishi space is the base of a deformation family $\m{M} \to V$ centred at $J$.

\begin{definition}
We say that the deformation family is \emph{complete} if, for another deformation family $\m{M}_1 \to V_1$, there exists a map
\[
\tau : V_1 \to V
\]
such that $\m{M}_1 = \tau^*\m{M}$. Moreover, if the differential of $\tau$ at the origin is unique we say that the deformation family is \emph{versal}, and \emph{universal} if $\tau$ itself is unique.
\end{definition}
We remark that Kuranishi's theorem \cite{Kuranishi_family_cpx_str} also gives the existence of a versal deformation family centred at $J$ which is complete for nearby complex structures.

\subsection{Fibrations}\label{Subsec:fibrations}
In this section, we collect some results on K\"ahler metrics on holomorphic submersions.
Let $\pi_X:X \to B$ be a holomorphic submersion of compact complex manifolds. Let $n$ be the dimension of $B$, $m$ the dimension of the fibres, so $\mrm{dim}(X)=n+m$.
We also assume that the spaces $H^0(X_b, T^{1,0}X_b)$ of holomorphic vector fields on $X_b$ are isomorphic as Lie algebras for all $b$.

\begin{definition}
A line bundle $H_X$ on the total space $X$ of the submersion $\pi_X$ is said to be \emph{relatively ample} if its restrictions to the fibres of $\pi_X$ are ample.
Similarly, a representative $\omega \in c_1(H_X)$ is called a \emph{relatively K\"ahler metric} if its restrictions to the fibres of $\pi_X$ are K\"ahler metrics.
\end{definition}

Let $\omega$ be a relatively K\"ahler metric on $X$ and $\omega_B$ a K\"ahler metric on $B$.
We can define a K\"ahler metric on $X$ by taking the relative K\"ahler metric $\omega$ and adding a large multiple of the pull-back of the base metric:
\begin{equation*}
\omega_k := \omega + k\omega_B \qquad k \gg0.
\end{equation*}
The relative symplectic form $\omega$ determines a splitting of the tangent space
\begin{equation}\label{Eq:splitting_tangent_bundle}
TX = \m{V} \oplus \m{H}^{\omega},
\end{equation}
where $\m{V}_x = T_x X_{\pi(x)}$ is the tangent space to the fibre $X_x$, and
\begin{equation*}
\m{H}^{\omega}_x = \set*{\left. u \in T_xX\right| \omega(u, v)=0 \ \forall v \in \m{V}_x}.
\end{equation*}
In the context of symplectic fibrations $\omega$ is called a \emph{symplectic connection} \cite[Chapter 6]{McDuff_Salamon} and the splitting \eqref{Eq:splitting_tangent_bundle} induces a splitting on all tensor bundles.

In local computations, we make the following notation conventions:
\begin{enumerate}[label=(\roman*)]
\item[] $\{ w^1, \dots, w^m \}$ denote vertical holomorphic coordinates; indices are denoted with the letters $a, b, c, \dots$;
\item[] $\{z^1, \dots, z^n\}$ denote holomorphic coordinates on the base; indices are denoted with the letters $i, j, k, \dots$;
\item[] $\{\zeta^1, \dots, \zeta^{n+m}\}$ denote holomorphic coordinates on the total space; indices are denoted with the letters $p,q,r, \dots$.
\end{enumerate}

\begin{lemma}\label{Lemma:contraction_1-1_tensor}
Let $\alpha$ be a covariant $2$-tensor of type $(1,1)$. Then
\begin{equation*}
\Lambda_{\omega_k} \alpha = \Lambda_{\m{V}}\alpha + \frac{1}{k}\Lambda_{\omega_B}\alpha + O\left(k^{-2}\right),
\end{equation*}
where $\Lambda_{\omega_k}$ denotes the contraction with the K\"ahler metric $\omega_k$.
\end{lemma}
\begin{proof}
At each point $x \in X_b \subset X$, the matrix $[\omega_k]$ is block-diagonal:
\begin{equation*}
\begin{pmatrix}
[\omega_b] &0\\
0 & [\omega_{k, \m{H}}]
\end{pmatrix}.
\end{equation*}
Therefore in local coordinates around the point $x$
\begin{equation*}
(\omega_k)^{p\bar{q}}\alpha_{p\bar{q}} = (\omega_b)^{a\bar{b}}\alpha_{a\bar{b}} + (\omega_{k,\m{H}})^{i\bar{j}}\alpha_{i\bar{j}}.
\end{equation*}
The horizontal part of $\omega_k$, denoted $\omega_{k, \m{H}}$ splits as $\omega_{\m{H}}+ k\omega_B$, where $\omega_{\m{H}}$ is the horizontal part of $\omega$.
Let $[\omega_{k, \m{H}}]$, $[\omega_{\m{H}}]$ and $[\omega_B]$ be the matrices of coefficients of the two-forms $\omega_{k, \m{H}}$, $\omega_{\m{H}}$ and $\omega_B$ respectively. We can write
\begin{equation*}
[\omega_{k,\m{H}}] = k \left( \frac{[\omega_{\m{H}]} [\omega_B]^{-1}}{k} + \idmatrix\right) [\omega_B],
\end{equation*}
where $\idmatrix$ is the identity matrix and the base form $\omega_B$ induces a Riemannian metric on the horizontal tangent bundle, so its inverse is well defined. The inverse of the matrix $[\omega_{k,\m{H}}]$ can be expanded in inverse powers of $k$ as
\begin{equation*}
\begin{aligned}
[\omega_{k,\m{H}}] ^ {-1} &= k^{-1}[\omega_B]^{-1}\left( \frac{[\omega_{\m{H}}] [\omega_B]^{-1}}{k} + \idmatrix\right)^{-1} \\
&=k^{-1}[\omega_B]^{-1} \sum_{i=0}^\infty \left( -\frac{[\omega_{\m{H}}] [\omega_B]^{-1}}{k}\right)^i \\
&= \frac{1}{k} [\omega_B]^{-1} + O\left( k^{-2} \right).
\end{aligned}
\end{equation*}
This implies the claim.
\end{proof}

Let us now restrict to the case when the fibres of $X$ each admit a cscK metric. We also assume that the dimension of the vector space of holomorphic vector fields on a fibre $H^0(X_b, T^{1,0}X_b)$ is independent of $b$.
The following lemma explains that one can define a relative K\"ahler metric on the total space which is relatively cscK.
\begin{lemma}[{\cite[Lemma 3.8]{DervanSektnan_OSC1}}]
For any $b \in B$, let $\omega_b$ be a cscK metric on the fibre $X_b$ in the class $c_1(H_X|_b)$. Then there exists a $\omega \in c_1(H_X)$ which is relatively cscK.
\end{lemma}

Under the relatively cscK assumption, we now explain how the function space $C^\infty(X, \bb{R})$ splits in a way that takes into account the fibration structure.
We begin by considering the vertical Lichnerowicz operator,
\begin{equation*}
\m{D}_{\m{V}}^*\m{D}_{\m{V}} : C^\infty(X, \bb{R}) \to C^\infty(X, \bb{R}),
\end{equation*}
defined fibrewise as $\left( \m{D}_{\m{V}}^*\m{D}_{\m{V}} \varphi \right)|_{X_b} = \m{D}_b^*\m{D}_b \left. \varphi\right|_{X_b}$. It is a real operator since the fibrewise metric is cscK. By integrating a function $\varphi \in C^\infty(X, \bb{R})$ over the fibres of $\pi$, we define a projection
\begin{equation*}
\begin{aligned}
C^\infty(X, \bb{R}) &\longrightarrow C^\infty(B, \bb{R}) \\
\varphi &\longmapsto \int_{X/B} \varphi \omega^m.
\end{aligned}
\end{equation*}
Its kernel is given by the space $C^\infty_0(X, \bb{R})$ of functions that have fibrewise mean value zero. A key step in the study of optimal symplectic connections is that we can further split this space as follows.

Consider the real vector bundle $E \to B$ \cite[\S3.1]{DervanSektnan_OSC1}, whose fibre over $b \in B$ is the real finite-dimensional vector space $\mrm{ker}_0(\m{D}_b^*\m{D}_b)$ of holomorphy potentials on the fibre $X_b$ with mean-value zero with respect to $\omega_b$. It is well defined as a vector bundle since we assume that the complex dimension of the space $H^0(X_b, T^{1,0}X_b)$ of holomorphic vector fields on $X_b$ is independent of $b$ \cite[\S2.3]{Hallam_geodesics}.
The space of smooth global sections of $E$, denoted by $C^\infty(E)$, is given by the kernel over the fibrewise mean-value-zero functions of the vertical Lichnerowicz operator $\m{D}_{\m{V}}^*\m{D}_{\m{V}}$.

\begin{remark}
In \cite[Lemma 2.7]{Hallam_geodesics}, Hallam used the Cartan decomposition for the space $\f{h}(X_b)$ of holomorphic vector fields of the fibre to show that $E_b$ can be also viewed as the vector space of all K\"ahler potentials $\varphi_b$ on $X_b$ of mean-value zero for which $\omega_b + i\del\delbar\varphi_b$ is still cscK.
\end{remark}

We can split $C^\infty_0(X)$ as
\begin{equation*}
C^\infty_0(X, \bb{R}) = C^\infty(E) \oplus C^\infty(R),
\end{equation*}
where $C^\infty(R)$ is the fibrewise $L^2$-orthogonal complement with respect to the fibre metric $\omega_b$, i.e.\ for all $\varphi \in \mrm{ker}_0 \m{D}_b^*\m{D}_b$, $\psi \in C^\infty(R)$
\begin{equation*}
\langle \varphi , \psi \rangle_b := \int_{X_b} \varphi \psi \omega_b^m = 0.
\end{equation*}
So we obtain
\begin{equation}\label{Eq:splitting_function_space}
C^\infty(X, \bb{R}) = C^\infty(B) \oplus C^\infty(E)\oplus C^\infty(R).
\end{equation}
We denote by $p_E : C^\infty(X) \to C^\infty(E)$ the projection.
\begin{definition}
We denote by $\m{K}_E$ the space of functions $\varphi \in C^\infty(X)$ such that $\omega + i \del\delbar\varphi$ is still a \emph{fibrewise cscK metric}.
\end{definition}
In particular, if we change the relatively cscK metric $\omega$ to $\omega + i\del\delbar\varphi$ with $\varphi \in \m{K}_E$, the vector bundles $E(\omega)$ and $E(\omega+i\del\delbar\varphi)$ are isomorphic.
The following result \cite[Lemma 4.20]{DervanSektnan_OSC3} relates the space $\m{K}_E$ to the vector bundle $E\to B$.
\begin{lemma}\label{Lemma:derivative_K_E}
Let $\varphi_t : [0, 1] \to \m{K}_E$ be a smooth path. Then for all $t$
\begin{equation*}
\dot{\varphi}_t \in C^\infty(B) \oplus C^\infty(E(\omega + i\del\delbar\varphi_t))
\end{equation*}
that is to say that, up to a function on the base, $\dot{\varphi_t}$ is a fibrewise holomorphy potential with respect to $\omega+i\del\delbar\varphi_t$.
\end{lemma}

\subsection{Optimal symplectic connections}\label{Subsec:OSC}
Let $(Y, H_Y)\to (B,L)$ be a holomorphic submersion with a relatively ample line bundle $H_Y \to Y$.
In this section we recall the theory of optimal symplectic connections \cite{Ortu_OSCdeformations}.
The following assumptions restrict the class of admissible fibrations to those whose fibres satisfy a stability property defined in terms of K-stability.
More precisely we assume that:
\begin{enumerate}[label = (\roman*)]
\item\label{Assumption1} the fibres $Y_b$ are analytically K-semistable, which means by definition that they each admit a degeneration to a cscK manifold $X_b$. We also assume that the deformation is compatible with the fibration structure in the following sense: there exists a holomorphic map $\varpi:(\m{X}, \m{H}) \to (B, L) \times S$, parametrised by a disk $S$, such that for $s\ne 0$, the family $(\m{X}_s, \m{H}_s) \to B$ is isomorphic to the original fibration $\pi_Y:(Y, H_Y) \to B$ and the central fibration at $s=0$ is a family $\pi_X:(X, H_X) \to B$ whose fibres are cscK;
\item\label{Assumption2} the automorphism groups $\mrm{Aut}_0(X_b, H_b)$ of the fibres are all isomorphic.
\end{enumerate}
The first hypothesis is a stability assumption. We will refer to the submersion $X \to B$ as the \emph{relatively cscK degeneration} of $Y\to B$.
The second hypothesis holds if and only if the spaces $H^0(X_b, T^{1,0}X_b)$ are isomorphic as Lie algebras, which we assumed in \S\ref{Subsec:fibrations} to define the vector bundle $E\to B$.
It is needed to define optimal symplectic connection and it is a key hypothesis for the construction of the moduli space.

A relative version of Ehresmann's theorem \cite[Proposition 4.5]{Ortu_OSCdeformations} implies that $X$ and $Y$ are diffeomorphic. Let $M$ denote the underlying smooth manifold.
Since the Chern classes are integral classes we have that $c_1(H_X)$ coincides with $c_1(H_Y)$ as cohomology classes on $M$. Since $Y$ is a small deformation of $X$, the cohomology class $c_1(H_X)$ is of type $(1,1)$ on $Y$, so $\omega$ is a $(1,1)$-form with respect to the complex structure of $Y$ \cite[\S6.1]{Huybrechts_ComplexGeometry}.
By Moser's theorem \cite[Theorem 7.2]{DaSilva_Lectures_symplectic_geometry}, we can modify the complex structure of $Y$ by a small diffeomorphism so that $\omega$ restricted to the fibres of $\pi_Y$ is compatible with the restriction of the complex structure. Thus we can assume that $\omega$ is relatively K\"ahler on $Y$.

Therefore we can view $Y\to B$ and $X\to B$ as the same relatively symplectic fibration $(M, \omega) \to B$ with two different integrable almost complex structures $J$ and $I$
where $(\omega, I)$ is relatively cscK and $(\omega, J)$ is just relatively K\"ahler. The family $\m{X} \to B\times S$ corresponds to a family of complex structures $\{J_s\}$ on $(M, \omega) \to B$, such that for $s \ne 0$, $J_s$ is isomorphic to $J$ and $J_0$ is isomorphic to $I$.

The definition of optimal symplectic connections involves various curvature quantities and a description of the complex structure of $Y$ as a deformation of the one on $X$.
To describe the deformations of the complex structure $I$, we begin by considering the space $\scr{J}_\pi$ of almost complex structures compatible with $\omega$ and such that $\dd\pi \circ J = J_B \circ \dd\pi$.
The tangent space at $I$ to $\scr{J}_\pi$ can be identified with
\[
T_{I}^{0,1}\scr{J}_\pi = \set*{\left.\alpha \in \Omega^{0,1}(\m{V}^{1,0}) \right| \omega_F (\alpha\cdot, \cdot) + \omega_F (\cdot, \alpha\cdot) = 0},
\]
where $\omega_F$ is the vertical part of $\omega$.
In local holomorphic coordinates an element $\alpha \in T_{J}^{0,1}\scr{J}_\pi$ is written as
\[
\alpha = \alpha^a_{\ \bar{b}} \del_{w^a}\otimes \dd \bar{w}^b + \alpha^a_{\ \bar{j}}\del_{w^a}\otimes \dd \bar{z}^j,
\]
where \cite[Lemma 4.1]{Ortu_OSCdeformations}
\begin{equation}\label{Eq:local_coords_deformations}
\alpha^{a}_{\ \bar{j}} = \alpha^a_{\ \bar{c}}(\omega_{F})^{d\bar{c}} (\omega)_{d\bar{j}}.
\end{equation}
Consider the map
\begin{equation*}
\begin{aligned}
P_{\m{V}}:C^\infty_0(X, \bb{R}) &\longrightarrow T_{I}^{0,1}\scr{J}_\pi \\
\varphi &\longmapsto \delbar (\mrm{grad}^{\omega_F}\varphi)^{1,0},
\end{aligned}
\end{equation*}
which is the relative version of the map \eqref{Eq:map_P}.
Let $\widetilde{H}^1_{\m{V}}$ be the kernel of the elliptic operator
\[
\square_{\m{V}} = P_{\m{V}}P_{\m{V}}^* + (\delbar^*\delbar)^2,
\]
where the adjoint is computed with respect to any K\"ahler metric on $X$ which restricts to $\omega_F$ vertically.  The operator $\square_{\m{V}}$ is fibrewise elliptic because $P_{\m{V}}P_{\m{V}}^*$ is elliptic in the fibre direction, and globally elliptic because $(\delbar^*\delbar)^2$ is.
The space $\widetilde{H}^1_{\m{V}}$ is the space of integrable first-order deformations of $I$. Let $K_\pi$ be the group of biholomorphisms of $I$ which restrict to an isometry on each fibre, with respect to the fibrewise metric defined by $(\omega, I)$. The group $K_\pi$ acts on $\widetilde{H}^1_{\m{V}}$ by pull-back.
The following theorem is a fibrewise version of Theorem \ref{Thm:Kuranishi}.

\begin{theorem}[{\cite[\S4.2]{Ortu_OSCdeformations}}]\label{Thm:relative_Kuranishi}
There exists a neighborhood of the origin $V_\pi \subset \widetilde{H}^1_{\m{V}}$ and a $K_\pi$-equivariant map
\begin{equation}\label{Eq:rel_Kuranishi_map}
\Phi : V_\pi \rightarrow \scr{J}_\pi
\end{equation}
such that:
\begin{enumerate}
\item $\Phi(0) = I$;
\item If $v_1, v_2 \in V_\pi$ and $v_1 \in K_b^{\bb{C}}\cdot  v_2$, and if $\Phi (v_1)$ is integrable, then $\Phi(v_1)$ and $\Phi(v_2)$ are isomorphic and for all $b \in B$ the complex structure $\left. \Phi(v_1) \right|_{X_b}$ is in the same $\scr{G}_b^c$-orbit as $\left. \Phi(v_2) \right|_{X_b}$;
\item For any $J \in \scr{J}_\pi$ integrable close to $I$, there exists $J'$ in the image of $\Phi$ such that, for all $b$, $J_b'$ is in the same $\scr{G}_b^c$-orbit as $J_b$;
\item $\mrm{Scal}_{\m{V}}\left(\omega_X, \Phi(x)\right) \in C^\infty(E, I)$.
\end{enumerate}
\end{theorem}

We now describe the curvature quantities involved in the definition of optimal symplectic connections.
The symplectic connection $\omega$ determines the following curvature quantities:
\begin{enumerate}[label=(\roman*)]
\item the \emph{symplectic curvature} is a two-form on $B$ with values in the fibrewise Hamiltonian vector fields defined for $v_1, v_2 \in \f{X}(B)$ as
\begin{equation*}
F_{\m{H}}(u_1, u_2) = [u_1^\sharp, u_2^\sharp]^{\mrm{vert}},
\end{equation*}
where $u_j^\sharp$ denotes the horizontal lift. Let $\gamma$ be the map which associates to a fibrewise Hamiltonian vector field its fibrewise Hamiltonian function with fibrewise mean value zero.
Thus we consider $\gamma^*(F_{\m{H}})$, which is a two-form on $B$ with values in $C^\infty_0(Y, \bb{R})$, and we pull it back to $Y$;
\item the curvature $\rho$ of the Hermitian connection induced by $(\omega, I)$ on the top wedge power $\ext{m}\m{V}$. We will primarily consider its purely horizontal part $\rho_{\m{H}}$;
\item the curvature of the deformation family $\nu$, which is defined as follows. On each fibre $X_b$ of the relatively cscK degeneration $X \to B$, let $\widetilde{H}^1_b$ be the deformation space \eqref{Eq:Htilde_1_mfd} and $V_b$ the Kuranishi space described in Theorem \ref{Thm:Kuranishi}. Let $v_b \in \widetilde{H}^1_b$ and denote by $\mbf{v}_b$ a vector field on $V_b$ such that $\mbf{v}_b\vert_{0} = v_b$. For any $f \in \f{k}_b$ let
\[
A_f v_b := -(\m{L}_{\mrm{grad}^\omega f} \mbf{v}_b)\vert_0.
\]
Define $\nu_b$ as
\begin{equation*}
\langle \nu_b(v_b), f \rangle = \frac{1}{2}\Omega_0(A_f v_b, v_b).
\end{equation*}
In particular, $\nu_b$ is the moment map for the action of the group $K_b$ on the tangent space at the origin to $V_b$.
Therefore we can define a global section $\nu$ of the vector bundle $E \to B$ of relatively holomorphy potentials of $X \to B$, as
\begin{equation}\label{Eq:map_nu_relative}
\nu(b) = \nu_b(v_b).
\end{equation}
Given $x \in V_\pi$, $v \in \widetilde{H}^1_{\m{V}}$, we will write $\nu_x(v)$
when we want to underline the dependence of the map $\nu$ on the complex structure $\Phi(x)$ and on the deformation $v$.
\end{enumerate}
\begin{definition}
\cite[\S3.3]{Ortu_OSCdeformations}
A relatively K\"ahler metric $\omega$ on $Y \to B$ is called an \emph{optimal symplectic connection} if
\begin{equation}\label{Eq:OSC}
p_E \left(\Delta_{\m{V}}(\Lambda_{\omega_B} m^*(F_{\m{H}}))+ \Lambda_{\omega_B}\rho_{\m{H}}\right) + \frac{\lambda}{2}\nu=0,
\end{equation}
for a positive number $\lambda$.
In the following, we will use the notation $\Theta(\omega, J) = \Delta_{\m{V}}(\Lambda_{\omega_B} m^*(F_{\m{H}}))+ \Lambda_{\omega_B}\rho_{\m{H}}$.
\end{definition}

Equation \eqref{Eq:OSC} is a second-order elliptic equation on the vector bundle $E \to B$.
The equation with $\nu=0$ is the condition for an optimal symplectic connection in the sense of \cite{DervanSektnan_OSC1}, where all the fibres are required to be cscK.

The linearisation of the equation at a solution is given by the operator \cite[\S5.3]{Ortu_OSCdeformations}
\[
\widehat{\m{L}} = \m{R}^*\m{R} + \m{A}^*\m{A}
\]
on $C^\infty(E)$, where
\begin{equation}\label{Eq:operator_R}
\m{R}(\varphi_E) = \delbar_B \nabla_{\m{V}}^{1,0} \varphi_E
\end{equation}
and
\begin{equation}\label{Eq:operator_mA}
\m{A}(\varphi_E) = \dd_0\Phi\left(A_{\varphi_E} v\right).
\end{equation}
The adjoint is computed with respect to $\omega_F + \omega_B$. Here $\nabla_{\m{V}}^{1,0} \varphi_E$ is a section of the holomorphic tangent bundle; the vertical part of $\delbar \nabla_{\m{V}}^{1,0} \varphi_E$ vanishes since $\varphi_E \in C^\infty(E)$ and the horizontal part is denoted by the expression \eqref{Eq:operator_R}.
The operator \eqref{Eq:operator_R} can be described as follows \cite[\S4.3]{DervanSektnan_OSC1}: let $\m{D}_k^*\m{D}_k$ be the Lichnerowicz operator with respect to the K\"ahler metric $\omega_k$. It admits a power series expansion in negative powers of $k$:
\begin{equation*}
\m{D}_k^*\m{D}_k = \m{L}_0 + k^{-1}\m{L}_1 + O\left(k^{-2}\right),
\end{equation*}
where $\m{L}_0$ is the \emph{vertical} Lichnerowicz operator $\m{D}_{\m{V}}^*\m{D}_{\m{V}}$. Then for $\varphi, \psi$ fibrewise holomorphy potentials
\begin{equation*}
\int_X \varphi \m{L}_1(\psi) \omega^m \wedge \omega_B^n = \int_X \langle \m{R}\varphi, \m{R}\psi \rangle_{\omega_F+\omega_B} \omega^m\wedge \omega_B^n.
\end{equation*}
This means that the operator $\m{R}^*\m{R}$ can actually be seen as $p_E \circ \m{L}_1$ restricted to $\m{C}^\infty_E(X)$. The kernel of $\m{R}$, thus of $\m{R}^*\m{R}$, consists of fibrewise holomorphy potentials which are global holomorphy potentials on $X$ with respect to $\omega_k$, and it is independent of $k$.
The operator \eqref{Eq:operator_mA} is described as
\begin{equation}\label{Eq:linearised_relative_map_nu}
\langle \dd_v \nu (\m{L}_{\nabla_{\m{V}}{\varphi}}v), \psi \rangle = \int_X \langle \dd_0 \Phi \left(A_{\varphi}v\right), \dd_0 \Phi\left(A_\psi v\right) \rangle_{\omega_F} \omega_F^m \wedge \omega_B^n,
\end{equation}
where $\varphi, \psi \in C^\infty(E)$. A function $\psi \in C^\infty(E)$ is in the kernel of $A$ if and only if $\psi$ is a fibrewise holomorphy potential with respect to all $J_s$, i.e. $\psi \in C^\infty(E, J_s)$.
Thus the kernel of $\widehat{\m{L}}$ consists of those functions $\psi \in C^\infty(E, J_0)$ such that  $\delbar_s (\nabla_{s, \m{V}}^{1,0}\psi )= 0$ for all $s$.


\section{Construction of the moduli space}\label{Sec:analytic_moduli_space}
Consider a submersion $Y\to B$ that degenerates to a relatively cscK holomorphic submersion $X\to B$, as described in \S \ref{Subsec:fibrations}.
We prove that the set of deformations of $Y \to B$ that still degenerate to a relatively cscK fibration (possibly different from $X\to B$) forms a locally closed analytic subset of the relative Kuranishi space.
We then prove that the solutions to the optimal symplectic connection equation \eqref{Eq:OSC} form an open set inside the locally closed subset of admissible deformations, which allow us to define a local moduli space of optimal symplectic connections.
Finally we glue the local moduli spaces and we prove that we obtain a global Hausdorff complex space which parametrises optimal symplectic connections.

\subsection{Openness of the setting}\label{Subsec:openness_setting}
Given a fibration $\pi_Y: (Y, H_Y) \to (B, L)$ with analytically K-semistable fibres, we have introduced in \ref{Assumption1} the assumption that there exists a degeneration of $\pi_Y$ to a fibration $\pi_X : (X, H_X) \to (B, L)$ such that the fibres of $\pi_X$ are cscK.
We have explained in \S \ref{Subsec:fibrations} that we can consider $Y\to B$ and $X\to B$ as the same symplectic fibration $\pi: (M, \omega) \to B$, and the degeneration as a deformation of a complex structure $J$ to $I$, where $(\omega, I)$ ha fibrewise constant scalar curvature.
The goal of this section is to understand for which deformations $J'$ of $J$ we can still find a relatively cscK degeneration $X'$, and to construct such a degeneration.

We begin by working locally in $B$.
Let $\m{U} \subseteq B$ be a coordinate open subset of $B$.
The fibre over the origin of $\m{U}$, denoted $M_0$, has a constant scalar curvature metric $(\omega_0, I_0)$.
Let $K$ be the group of Hamiltonian isometries of $(\omega_0, I_0)$ and let $\widetilde{H}^1_0$ be the vector space \eqref{Eq:Htilde_1_mfd} parametrising first-order deformations of $(M_0, \omega_0, I_0)$.
Since $(\omega_0, I_0)$ has constant scalar curvature, the complexification of $K$ is the group
\[
G=\Aut_0(M_0, H_0).
\]

We will employ the following definition of GIT-stability for the action of a reductive group on an affine space \cite[\S3]{Newstead_ModuliproblemsOrbitspaces} \cite[\S4.5, \S4.6]{Hoskins_notes2015}.

\begin{definition}
Let $\bb{A}^d$ be an affine space of dimension $d$ and $G$ a reductive affine group acting on it. Let $(z_1, \dots, z_d)$ be a system of coordinates on $\bb{A}^d$. The space $\bb{A}^d$ can be embedded in the projective space $\bb{P}^d$ as a coordinate chart with the map
\[
(z_1, \dots, z_d) \mapsto [1:z_1: \dotsc :z_d].
\]
We extend the action of $G$ to an action on $\bb{P}^d$ by acting trivially on the first coordinate.
Let $x \in \bb{A}^d$ and let $\widehat{x}$ be its image in $\bb{P}^d$.
We say that $x$ is
\begin{enumerate}[label = (\roman*)]
\item \emph{semistable} if there exists a non-constant homogeneous polynomial in $[z_0: \dotsc: z_d]$ on $\bb{P}^d$ which is $G$-invariant and does not vanish at $\widehat{x}$;
\item \emph{polystable} if it is semistable and its orbit is closed;
\item \emph{stable} if it is semistable, its orbit is closed and its stabiliser is finite.
\end{enumerate}
\end{definition}
In particular, the polynomial $P(\mbf{z}) = z_0$ is a $G$-invariant homogeneous polynomial that does not vanish at any point of $\bb{A}^d$. So every point of $\bb{A}^d$ is semistable.

\begin{remark}
The fact that every point is semistable holds because we defined the action on $\bb{P}^d$ to be trivial on the first homogenous coordinate. In principle, one can choose to extend the action in a non-trivial way and can still define stability as above, but unstable points might appear.
We do not treat this case here, so we have included the trivial extension of the action in the definition of stability.
\end{remark}

We apply the definition of stability to the vector space $\widetilde{H}^1_0$.
The complex structure $I_0$ corresponds to the origin in the Kuranishi space $V_0$, which is fixed by the action of $G$.
Therefore its orbit is closed and its stabiliser is the group $G$ itself, so it is a polystable point.
A key result for our construction is the following.
\begin{lemma}[{\cite[Corollary 4.31]{Hoskins_notes2015}}]\label{Lemma:unique_ps_orbit}
The closure of the orbit of every point in $\widetilde{H}^1_0$ contains a unique polystable orbit.
\end{lemma}

Let $V_0$ be the subspace of the Kuranishi space which parametrises integrable almost complex structures; it is a locally closed analytic subspace of $\widetilde{H}^1_0$ because it is defined by the vanishing of the Nijenhuis tensor.
Thus the family $X \to B$ can be described locally over $\m{U}$ as a family $\{ x_b \}$ in $V_0$.
By our hypothesis \ref{Assumption2}, the automorphism group of $(M_b, H_b)$ is isomorphic to $G$ for all $b \in B$. Therefore the points $\{ x_b \mid b \in \m{U} \}$ are all fixed by the action of $G$ and are hence polystable.

The family $Y \to B$ can be described locally over $\m{U}$ as a family $\{ y_b\mid b \in \m{U}  \}$ of points such that for each $b$ the closure of the $G$-orbit of $y_b$ contains the polystable point $x_b$ \cite[Theorem 1.3]{ChenSun_CalabiFlow}.
By Lemma \ref{Lemma:unique_ps_orbit}, $x_b$ is the only polystable orbit in the closure of the orbit of $y_b$.
Let $V_0^+$ be the set of all semistable points that have a fixed point in the closure of their orbit.
Then the map
\begin{equation}\label{Eq:ss-ps}
\begin{aligned}
F: V_0^+ &\to V_0
\end{aligned}
\end{equation}
that maps a semistable point to the corresponding fixed point is well-defined.

\begin{lemma}\label{Lemma:local_construction}
The set $V_0^+$ is an analytic subvariety of $V_0$ and the map \eqref{Eq:ss-ps} is holomorphic.
\end{lemma}
\begin{proof}
The space $V_0$ is an open subset of the vector space $\widetilde{H}^1$. Let $d$ be the dimension of $\widetilde{H}^1$, so each point $z \in V_0$ has coordinates
\[
\left(z_1, \dots, z_d\right).
\]

Let us begin with the case when $G$ is isomorphic to a complex torus ${(\bb{C}^*)}^r$.
The fixed points of the action can be described by the vanishing of the coordinates
\[
z_{i_1} = \dots = z_{i_h} = 0 \quad\text{for}\quad i_1, \dots, i_h \in \{1, \dots, d\},
\]
thus they form an analytic subspace of $V_0$.

By the Hilbert-Mumford criterion \cite[Theorem 2.1]{MumfordFogartyKirwan_GIT}, a point $y_j$ is semistable if it is semistable for the action of every 1-parameter subgroup of the group $G$.
Let $\rho(t): \bb{C^*} \hookrightarrow G$ be a 1-parameter subgroup. The action of $\rho(t)$ can be written as $(t^{a_1}, \dots, t^{a_d})$, where the numbers $a_j$ are the \emph{weights} of the action. Then $V_0$ splits into a sum of weight spaces $V_0^{\mrm{pos}} \oplus V_0^{\mrm{fix}} \oplus V_0^{\mrm{neg}}$, where $\rho(t)$ acts on $V_0^{\mrm{pos}}$ with positive weights, on $V_0^{\mrm{neg}}$ with negative weights and fixes the subspace $V_0^{\mrm{fix}}$.
A semistable point that has a fixed point in the closure of its orbit is described by the following condition: if a coordinate in $V_0^{\mrm{pos}}$ is nonzero, then all coordinates in $V_0^{\mrm{neg}}$ vanish.
Applying this condition to all 1-parameter subgroups yields a set of polynomial equations that define the semistable points in $V_0^+$. Thus, the semistable points correspond to an analytic subset of $V_0$.
The map $F$ is the projection onto the set described by $\{ z_{i_1} = \dots \ = z_{i_h} = 0\}$, thus it is holomorphic.

Let now $G$ be any reductive group, and let $y$ be a semistable point and $x$ be a polystable and fixed point in the closure of its orbit.
The fixed points of $G$ form a vector subspace $V_0^\mrm{fix}$ also in this case.
Moreover there exists a 1-parameter subgroup $\lambda_x : \bb{C}^* \hookrightarrow G$ such that
\[
\lim_{t\to 0} \lambda_x(t)\cdot y = x.
\]
It follows that the map $y \mapsto x$ is the projection onto the vector subspace of fixed points for $\lambda_x$, hence it is holomorphic.
Consider the composite map
\begin{equation}\label{Eq:composite_map_proj}
V_0 \to V_0 \overset{pr}{\to} V_0^{\mrm{fix}}
\end{equation}
where the first map is the projection onto the vector subspace of fixed points for $\lambda_x$ and the second map is the projection onto the subspace of fixed points for the whole group $G$.
The map is holomorphic because it is a composition of holomorphic projections. We prove that it coincides with the map \eqref{Eq:ss-ps}.
Let $\widetilde{x}' = \lim_{t\to 0} \lambda_x(t)\cdot y'$. Then
\[
\overline{G\cdot \widetilde{x}'} \subseteq \overline{G\cdot y'}.
\]
The unique polystable orbit contained in $\overline{G\cdot \widetilde{x}'}$ is also a polystable orbit in $\overline{G\cdot y'}$, so it must coincide with the fixed-point orbit $\{x'\}$.
Thus flowing along the orbit of $\widetilde{x}'$ amounts to projecting onto the subspace of fixed points, and so the map \eqref{Eq:composite_map_proj} maps any point $y'$ to the fixed point $x'$ in the closure of its orbit.
\end{proof}

\begin{remark}
The map \eqref{Eq:ss-ps} is analogous to the one given by the Bya\l inicki-Birula decomposition \cite{BB_decomposition,BB_decomposition_reductive_groups}.
\end{remark}


Let $V_\pi$ be the Kuranishi space of the fibration $X \to B$ defined in Theorem \ref{Thm:relative_Kuranishi}. Consider the subspace
\[
V^+_\pi := \set*{y \in V_\pi \mid y\vert_{X_b} \in V_b^+}.
\]
We remark that $V^+_\pi$ depends on the complex structure of the reference fibration $X\to B$ and its deformation $Y\to B$.
We denote by $\scr{J}^+_\pi$ the image of $V_\pi^+$ via the relative Kuranishi map \eqref{Eq:rel_Kuranishi_map}.

\begin{lemma}\label{Lemma:gluing_V+b}
$V_\pi^+$ is a locally closed subvariety of $V_\pi$.
\end{lemma}
\begin{proof}
Let $\m{U}_1$ and $\m{U}_2$ be two open coordinate subsets of $B$ with non-empty intersection and let $0_1$ be the origin of $\m{U}_1$ and $0_2$ be the origin of $\m{U}_2$.
Denote by $I_1$ and $I_2$ the complex structures of the fibres $X_{0_1}$ and $X_{0_2}$, and $G_1$, $G_2$ their groups of automorphisms. By assumption, $G_1$ and $G_2$ are isomorphic, and we will denote them by $G$.
Recall that the Kuranishi space is \emph{versal} and more specifically that it is a complete deformation space for the nearby fibres.
We use the versality of the Kuranishi map to glue the spaces $V_1^+$ and $V_2^+$ constructed in Lemma \ref{Lemma:local_construction} to a variety $V^+$ and prove that $V_\pi^+$ is obtained as the intersection of said variety with the relative Kuranishi space $V_\pi$.
More precisely, versality of the Kuranishi space means that there is a map
\[
\tau_{21} : V_2 \to V_1,
\]
not necessarily unique.

The map $\tau_{21}$ can be taken to be $G$-equivariant.
In fact, the $G$-equivariance can be traced back to the proof of Kuranishi's Theorem \ref{Thm:Kuranishi}. The map $\tau_{21}$ is defined using the implicit function theorem, which can be applied to a $K$-equivariant map to provide an implicit inverse function which is $K$-equivariant. Since $G$ is the complexification of $K$, we obtain that $\tau_{21}$ is $G$-equivariant.
The equivariance also implies that the image of $V_2^+$ is $V_1^+$, so we can restrict $\tau_{21}$ to
\[
\widetilde{\tau}_{21}:V_2^+ \to V_1^+.
\]
Moreover the map $\widetilde{\tau}_{21}$ has an inverse that is constructed reversing the roles of $V_1$ and $V_2$, so it is an isomorphism.
In fact, although the map $\tau_{21}$ is not canonical, the restriction to $\widetilde{\tau}_{21}$ is fixed by the reference K-semistable fibration $Y\to B$.
Each Kuranishi space $V_b$ is a complex subspace of the vector space $\widetilde{H}_b^1$, described as the kernel of the elliptic operator $P_bP_b^*+(\delbar_b^*\delbar_b)^2$ \eqref{Eq:Htilde_1_mfd}.
So we can use the isomorphism $\widetilde{\tau}_{21}$ to glue the spaces $V_b^+$ to a subvariety $V^+$ of the kernel of the fibrewise elliptic operator
\[
P_{\m{V}}P_{\m{V}}^*+ (\delbar^*_{\m{V}}\delbar_{\m{V}})^2.
\]
Therefore the intersection
\[
V_\pi^+ = V^+ \cap V_\pi
\]
is a locally closed subvariety of $V_\pi$.
\end{proof}

The following lemma shows that we can glue the local fibration constructed in Lemma \ref{Lemma:local_construction} to a global fibration over $B$.
\begin{lemma}
Let $Y'=(M, \omega, J')\to B$ be a fibration with complex structure $J'$ represented by $y'\in V^+_\pi$. Then $Y'$ degenerates to $X'=(M, \omega, I') \to B$ such that
\begin{enumerate}
\item $(\omega, I')$ is relatively cscK;
\item the groups $\mrm{Aut}(X'_b, H_b')$ are isomorphic for all $b \in B$.
\end{enumerate}
\end{lemma}
\begin{proof}
Let $\m{U} \subseteq B$ be an open coordinate subset. By the relative Kuranishi Theorem \ref{Thm:Kuranishi}, for each $b\in\m{U}$ there exists a point $y'_b \in V_0$ such that $\Phi_0(y'_b)$ is in the same $\m{G}^c$-orbit of $J'_b$. The fact that the map \eqref{Eq:ss-ps} is holomorphic implies that we can find polystable points $\{x'_b\}$ such that $\{\Phi_0(x'_b)\}$ are a holomorphic family of cscK complex structures over $\m{U}$ that are deformations of $\{J'_b\}$. Then we can construct a local relatively cscK fibration from the pullback diagram
\begin{equation}\label{Eq:universal_family_pullback}
\begin{tikzcd}
X_{\mathcal{U}}' \arrow[d] \arrow[r] & \mathcal{M}_{\mathcal{U}} \arrow[d, "p_{\mathcal{U}}"] \\
\mathcal{U} \arrow[r, "i"']          & V_0                                                   
\end{tikzcd}
\end{equation}
where $p_{\m{U}} : \m{M}_{\m{U}} \to V_0$ is Kuranishi's versal family and $i(b) = x_b'$.

Now we glue the local fibrations to a fibration $X'\to B$. Let $\m{U}_1$ and $\m{U}_2$ and $G$ as in Lemma \ref{Lemma:gluing_V+b}.
The associated Kuranishi maps are denoted respectively by $\Phi_1:V_1\to \scr{J}$ and $\Phi_2:V_2\to\scr{J}$.

For $b \in \m{U}_1 \cap \m{U}_2$, consider the complex structure $J_b'$. Since $J_b'$ can be regarded as a deformation of both $I_1$ and $I_2$, there exist points $y'_{b,1} \in V_1$ and $y'_{b,2} \in V_2$ such that
\[
\Phi_1(y'_{b,1}) = J_b' = \Phi_2( y'_{b,2} ).
\]
The diagram \eqref{Eq:universal_family_pullback} produces two fibrations $X_1' \to \m{U}_1$ and $X_2' \to \m{U}_2$.
In order to glue the local fibrations we need to prove that there is an isomorphism
\[
\phi_{12}:X'_1\vert_{\m{U}_1 \cap \m{U}_2} \overset{\sim}{\rightarrow} X'_2\vert_{\m{U}_1 \cap \m{U}_2}
\]
and that it satisfies the cocycle condition
\[
\phi_{23} \circ\phi_{12}=\phi_{13}
\]
on a triple intersection $\m{U}_1\cap\m{U}_2\cap\m{U}_3$ of open subsets of $B$.
In particular, the first condition produces a global compact complex manifold $X'$, while the cocycle condition implies that $X'$ admits a submersion onto $B$.

To prove the existence of the isomorphism $\phi_{12}$ we use again the fact that the Kuranishi space induces complete deformations on the nearby fibres. The map
\[
\tau_{21} : V_2 \to V_1
\]
is such that $\m{M}_2 = \tau_{21}^*\m{M}_1$. In particular, we have the following diagram
\[
\begin{tikzcd}
\mathcal{M}_2 \arrow[d]     &                                                                                 & \mathcal{M}_1 \arrow[d] \\
V_2 \arrow[rr, "\tau_{21}"] &                                                                                 & V_1                     \\
X_2' \arrow[r]              & \mathcal{U}_1\cap\mathcal{U}_2 \arrow[lu, "i_2", hook] \arrow[ru, "i_1"', hook] & X_1'. \arrow[l]         
\end{tikzcd}
\]
If we show that the diagram is commutative, i.e. $\tau_{21}\circ i_2 = i_1$, then
\[
X_2'\vert_{\m{U}_1\cap\m{U}_2} = i_2^*\m{M}_2 \simeq (\tau_{21}\circ i_2 )^* \m{M}_1 = i_1^*\m{M}_1 = X_1'\vert_{\m{U}_1\cap\m{U}_2}.
\]
The commutativity follows from the fact that the map $\tau_{21}$ is $G$-equivariant. Indeed the points $i_1(b)$ and $i_2(b)$ are defined as the fixed-point limits of semistable orbits. Moreover, an equivariant map between two spaces on which there is an action of the same group $G$ sends fixed points to fixed points and the closure of the orbit of $i_1(b)$ to the closure of the orbit of $i_2(b)$.

We now show the cocycle condition. Let $\m{U}$ be a triple intersection $\mathcal{U}_1\cap\mathcal{U}_2\cap\mathcal{U}_3$ and consider the following diagram
\[
\begin{tikzcd}
V_1 & V_2 \arrow[l, "\tau_{21}"]                                                                                              & V_3 \arrow[l, "\tau_{32}"] \arrow[ll, "\tau_{31}"', bend right] \\
    & \m{U} \arrow[lu, "i_1", hook] \arrow[ru, "i_3"', hook] \arrow[u, "i_2", hook] &                                                                
\end{tikzcd}.
\]
On $V_3$ we have two families pulled-back from $V_1$, namely $\tau_{31}^* \m{M}_1$ and $ (\tau_{21}\circ\tau_{32})^*\m{M}_1$.
They induce two distinct families on $\m{U}$, pulled-back using $i_3$.
Although in general it is \emph{not} true that $\tau_{31}$ is equal to the composition $\tau_{21}\circ\tau_{32}$, the commutativity of the arrows proved above implies that
\[
\tau_{21}\circ\tau_{32}\circ i_3 = \tau_{31}\circ i_3.
\]
Therefore the two families $\tau_{31}^* \m{M}_1$ and $ (\tau_{21}\circ\tau_{32})^*\m{M}_1$ coincide.
\end{proof}

\subsection{Openness of the space of optimal symplectic connections}\label{Subsec:openness_OSC}
Let $Y \to B$ be a fibration with K-semistable fibres, and assume that it degenerates to a fibration $X \to B$ with cscK fibres, in the sense of \S \ref{Subsec:fibrations}.
Let $Y'\to B$ be a deformation of $Y \to B$ in $\scr{J}^+_\pi$. Then $Y'\to B$ admits a degeneration to $X' \to B$, whose fibres are cscK, as explained in \S \ref{Subsec:openness_setting}.
The goal of this section is to show that if $Y$ admits an optimal symplectic connection then $Y'$ also does.

We denote by $I$ the complex structure of $X$, by $J$ the complex structure of $Y$ and we assume that $Y \to B$ is generated by $v_0 \in V_\pi$.
We also assume that $(\omega, J)$ is an optimal symplectic connection. 
Let $V_\pi^+$ the subvariety of $V_\pi$ which describes the family of complex structures $\scr{J}_\pi^+$.

\begin{proposition}
For every $\varphi \in \m{K}_E(I)$ there exists $f \in \mrm{Diff}_0(M)$ such that $f^*\omega_\varphi = \omega$ and $(M, \omega_\varphi, I)\to B$ is isomorphic to $(M, \omega, f^*I)\to B$.
\end{proposition}
\begin{proof}
Let us consider a potential $\varphi \in \m{K}_E(I)$ and a path $\{ \varphi_t\}$ in $\m{K}_E(I)$ from 0 to $\varphi$. \cite[Theorem 3.3]{Hallam_geodesics} guarantees that this path exists and that it is smooth. We define the relatively cscK metrics
\begin{equation*}
\omega_{t} = \omega + 2i\del\delbar\varphi_t
\end{equation*}
and the K\"ahler metrics
\begin{equation*}
\omega_{k,t} = \omega_{t}+ k\omega_B,
\end{equation*}
where the $\del, \delbar$ operators are taken with respect to the relatively cscK complex structure $I$. From Proposition \ref{Lemma:derivative_K_E}, we have that $\dot{\varphi}_t \in C^\infty(E(\omega_{t}, I))\oplus C^\infty(B)$, for all $t$. Thus the fibrewise Hamiltonian vector fields
\begin{equation*}
\eta_{t} :=\mrm{grad}^{\omega_{t}}{\dot{\varphi}_t}
\end{equation*}
are well-defined. Consider the vertical vector fields
\begin{equation*}
\xi_t := (\nabla^{g_t}\dot{\varphi}_t)_{\m{V}} = (I\mrm{grad}^{\omega_{t}}{\dot{\varphi}_t})_{\m{V}}.
\end{equation*}
Then fibrewise
\begin{equation}
\frac{\dd}{\dd t} \omega_{t} = -\m{L}_{I\eta_{t}}\omega_{t} = \m{L}_{\xi_t}\omega_{t}.
\end{equation}
Let $\{f_t, t\in [0,1]\}$ be the \emph{isotopy} of the time-dependent vector field $\xi_t$, i.e.\ the collection of diffeomorphisms of $M$ such that
\begin{equation*}
\frac{\dd}{\dd t} f_t = \xi_t(f_t), \quad f_0= \mrm{id}.
\end{equation*}
Since $\xi_t$ is vertical, $f_t \in \mrm{Diffeo}(M, \pi)$. Recall the following property \cite[6.4]{DaSilva_Lectures_symplectic_geometry}, for a smooth family $\sigma_t$ of tensors:
\begin{equation*}
\frac{\dd}{\dd t} f_t^*\sigma_t = f_t^* \left( \m{L}_{\xi_t} \sigma_t + \frac{\dd\sigma_t}{\dd t}\right),
\end{equation*}
where $\xi_t$ and $f_t$ are a time-dependent vector field and its isotopy respectively. Applying it to $\omega_{t}$ gives the fibrewise relation
\begin{equation*}
\frac{\dd}{\dd t} f_t^*\omega_{t} = 0,
\end{equation*}
which implies
\begin{equation*}
f_t^* \omega_{t} = f_0^*\omega = \omega.
\end{equation*}
Define $J_t = f_t^*I$. Then, for $t=1$, the two metrics $g(\omega_X, f_1^*I)$ and $g(\omega_X+2i\del\delbar\varphi, I)$ are fibrewise isometric, i.e.
\begin{equation*}
(M, \omega, f_1^*I) \simeq (M, \omega+ 2i\del\delbar\varphi, I)
\end{equation*}
as relative K\"ahler manifold with fibrewise constant scalar curvature.
\end{proof}

Thus, given an integrable relatively cscK complex structure $I$ in $\scr{J}_\pi$, we have a map
\begin{equation}\label{Eq:map_F_pi}
\begin{aligned}
F_\pi : \m{K}_E(I) &\longrightarrow \scr{J}_\pi \\
\varphi &\longmapsto f_1^*I =: F_\pi(\varphi, I)
\end{aligned}
\end{equation}
which locally parametrises all integrable complex structures in the same diffeomorphism class of $I$ that are fibrewise cscK with respect to the fixed $\omega$. Its differential at the origin is
\begin{equation*}
\dd_0F_\pi(\dot{\varphi}) = \left.\frac{\dd}{\dd t}\right|_{t=0} J_t = \left.\frac{\dd}{\dd t}\right|_{t=0}  f_t^*I = \left.\frac{\dd}{\dd t}\right|_{t=0}\m{L}_{\xi_t} I = \m{L}_{(I\mrm{grad}^{\omega}{\dot{\varphi}})_{\m{V}}} I = -\frac{1}{2} \delbar (\mrm{grad}^{\omega}{\dot{\varphi}})_{\m{V}},
\end{equation*}
where, again by Proposition \ref{Lemma:derivative_K_E}, $\dot{\varphi} \in C^\infty(E(\omega, I)) \oplus C^\infty(B)$.

Now let $v_0 \in V_\pi$ be the deformation of the complex structure which generates the family $Y \to B$ and let $\Phi$ be the relative Kuranishi map \eqref{Eq:rel_Kuranishi_map}. Then we can extend the map \eqref{Eq:map_F_pi} to a map
\[
\begin{aligned}
F'_\pi : \m{K}_E(I) &\to \widetilde{H}^1_{\m{V}} \\
\varphi &\mapsto f_1^*v_0 =: F_\pi' (\varphi, I, v_0).
\end{aligned}
\]
Its differential at the origin computed at $\dot{\varphi} \in C^\infty(E(\omega, I)) \oplus C^\infty(B)$ is
\begin{equation}\label{Eq:derivative_F_pi'}
\dd_0F'_\pi (\dot{\varphi}) = \left.\frac{\dd}{\dd t}\right\vert_{t=0} f_t^*v_0 = \m{L}_{(I\mrm{grad}^{\omega}{\dot{\varphi}})_{\m{V}}} v_0.
\end{equation}

\begin{definition}\label{Def:HandG}
We denote with
$\m{P}$ the set of triples $(\varphi, x, v) \in C^\infty(M)\times TV_\pi$ such that $\varphi \in \m{K}_E(\Phi(x))$ and $x$ is $K_\pi^{\bb{C}}$-polystable.
The optimal symplectic connection operator is the map
\[
\begin{aligned}
\m{G} : \m{P} &\to C^\infty(X) \\
(\varphi, x, v) &\mapsto p_{E(\varphi, x)} \left( \Theta(\omega, F_\pi(\varphi,\Phi(x)))\right) + \frac{\lambda}{2}\nu_{\varphi, x}\left(F'_\pi(\varphi, \Phi(x), v)\right).
\end{aligned}
\]
In this expression $E(\varphi, x)$ is the vector bundle of fibre holomorphy potentials with respect to the K\"ahler structure $(\omega, F_\pi(\varphi, \Phi(x)))$ and $\nu_{\varphi,x}$ is the map \eqref{Eq:map_nu_relative} computed with respect to the complex structure $F_\pi(\varphi,\Phi(x))$.
\end{definition}
We now compute the differential of $\m{G}$ along the $\varphi$-variable computed at $(0,0,v_0)$.
We write
\[
\m{G} (\varphi, x, v) = \m{G}_1(\varphi, x) + \m{G}_2(\varphi, x, v),
\]
and we split the computation into two separate lemmas.

\begin{lemma}\label{Lemma:Derivative_mG1}
Let $\varphi \in C^\infty(E(\omega, I)) \oplus C^\infty(B)$. The differential along the first variable of $\m{G}_1$ is
\[
D_1 \m{G}_1 \vert_{(0,0, v_0)}(\varphi) = -\frac{1}{2}\m{R}^*\m{R}(\varphi).
\]
\end{lemma}
\begin{proof}
Let $\{J_t\}$ be a family of relatively cscK complex structures compatible with $\omega$, such that $J_0 = I$ and consider the scalar curvature of $(\omega_k, J_t)$,
\begin{equation}\label{Eq:expansion_scalar_curvature}
\mrm{Scal}(\omega_k, J_t) = \mrm{Scal} (\omega_k, I) + t\left.\frac{\dd}{\dd t}\right|_{t=0} \mrm{Scal}(\omega_k, J_t) + O(t^2).
\end{equation}
Let $\alpha = \del_{t=0} J_t$, and define
\begin{equation*}
Q_k(\alpha) := \left.\frac{\dd}{\dd s}\right|_{s=0} \mrm{Scal}(\omega_k, J_t).
\end{equation*}
From Lemma \ref{Lemma:PofQ} we obtain
\begin{equation*}
\begin{aligned}
Q_k(\alpha) &= \mrm{Im} \left((g_k)^{p\bar{q}} \ \nabla_p \nabla_a \alpha^a_{\ \bar{q}}\right) = \mrm{Re} \left((\omega_k)^{p\bar{q}} \ \nabla_p \nabla_a \alpha^a_{\ \bar{q}}\right).
\end{aligned}
\end{equation*}
We need to compute the sub-leading order term of $Q_k$ along the differential of the map $F_\pi$ \eqref{Eq:map_F_pi}

Let us consider $\alpha = \delbar \left(\mrm{grad}^\omega {\varphi_E}\right)^{1,0}_{\m{V}}$, where $\varphi_E \in C^\infty_E(X)$ is a fibrewise holomorphy potential. Then
\begin{equation}\label{Eq:Q_Kspecial}
Q_k(\alpha) = k^{-1} \mrm{Re}\left( \m{R}^*\m{R} (\varphi_E\right)) + \Ok{-2},
\end{equation}
where $\m{R}(\varphi_E) = \delbar_B \mrm{grad}_{\m{V}} \varphi_E$ and the adjoint is computed with respect to $\omega_F + \omega_B$.
As explained in Section \ref{Subsec:OSC}, the operator $\m{R}^*\m{R}$ can actually be seen as $p_E \circ \m{L}_1$ restricted to $\m{C}^\infty(E(\omega,I))$. Its kernel consists of fibrewise holomorphy potentials which are global holomorphy potentials on $X$ with respect to $\omega_k$.

A local coordinate expression for $\alpha$ is
\begin{equation*}
\begin{aligned}
\alpha = \delbar \left(Y_{\varphi_E}\right)^{1,0}_{\m{V}} &= \del_{\bar{z}^j} \left(\omega_F^{a\bar{b}} \del_{\bar{w}^b}\varphi_E \right)\del_{w^a} \otimes \dd \bar{z}^j = \nabla_{\bar{j}} \left(\omega_F^{a\bar{b}} \nabla_{\bar{b}}\varphi_E \right)\del_{w^a} \otimes \dd \bar{z}^j .
\end{aligned}
\end{equation*}
In particular, since the potential $\varphi_E$ is a fibrewise holomorphy potential, the component of $\alpha$ with the covariant index of vertical type vanishes. Therefore, using Lemma \ref{Lemma:contraction_1-1_tensor},
\begin{equation*}
\begin{aligned}
Q_k(\alpha) 
&=  k^{-1}\mrm{Re} \left( \omega_B^{i\bar{j}} \ \nabla_i \nabla_a \nabla_{\bar{j}} \left(\omega_F^{a\bar{b}} \nabla_{\bar{b}}\varphi_E \right) \right) +  \Ok{-2} =\\
&=  k^{-1}\mrm{Re} \left( \omega_B^{i\bar{j}} \ \omega_F^{a\bar{b}} \ \nabla_i \nabla_a \nabla_{\bar{j}} \nabla_{\bar{b}}\varphi_E  \right) + \Ok{-2} =\\
&=k^{-1} \mrm{Re}\left( \m{R}^*\m{R} (\varphi_E\right)) + \Ok{-2}.
\end{aligned}
\end{equation*}

A similar computation also holds if we consider, instead of a potential in $C^\infty(E(\omega, I))$, a map $\varphi \in C^\infty(E(\omega, I))\oplus C^\infty(B)$, so that $\alpha=\delbar(Y_\varphi)^{1,0}_{\m{V}}$. This choice amounts to considering an element in the image of the differential of the map $F_\pi$ \eqref{Eq:map_F_pi}. In this case
\begin{equation*}
Q_k(\alpha) = \mrm{Re}(\m{L}_0\varphi) + \frac{1}{k}\mrm{Re}(\m{L}_1\varphi) + \Ok{-2},
\end{equation*}
where $\m{L}_0(\varphi)=\m{L}_0(\varphi_B) =0$, since $\varphi_B$ is constant when restricted to a fibre, so we obtain the same equation as \eqref{Eq:Q_Kspecial}.
To obtain the correct coefficient in the claimed expression note that 
\[
\dd_0F_\pi(\varphi) = -\frac{1}{2}\alpha. \qedhere
\]
\end{proof}

\begin{lemma}\label{Lemma:Derivative_mG2}
Let $\varphi \in C^\infty(E(\omega, I)) \oplus C^\infty(B)$. The differential along the first variable of $\m{G}_2$ is
\[
D_1 \m{G}_2 \vert_{(0,0, v_0)}(\varphi) = -\m{A}^*\m{A}(\varphi).
\]\end{lemma}
\begin{proof}
We compute
\[
\left.\frac{\dd}{\dd t}\right\vert_{t=0} \nu_t \left(f_t^*v_0\right),
\]
where $\nu_t$ is the map $\nu$ computed with respect to the complex structure $f_t^*\Phi(0)$.
Now, $f_t$ is the isotopy of the vector field $\xi_t = (I\mrm{grad}^{\omega_t}\dot{\varphi}_t)_{\m{V}}$, where $\dot{\varphi}_t$ is in $C^\infty(E_t)\oplus C^\infty(B)$ and $\dot{\varphi}_0 = \varphi$.
In the expression of $\xi_t$ we are fixing the complex structure $I$ and varying the K\"ahler form $\omega_t$.
In particular, $\xi_t$ is a fibrewise holomorphic vector field with respect to $I$.
This implies that $f_t^*I = I$, so $\nu_t = \nu$.
Therefore
\[
\left.\frac{\dd}{\dd t}\right\vert_{t=0} \nu_t \left(f_t^*v_0\right) = \dd_{v_0}\nu \left(\left.\frac{\dd}{\dd t}\right\vert_{t=0}f_t^*v_0 \right).
\]
Using the expression \eqref{Eq:derivative_F_pi'} we obtain
\[
\dd_{v_0}\nu \left(\left.\frac{\dd}{\dd t}\right\vert_{t=0}f_t^*v_0 \right) = \dd_{v_0}\nu \left( \m{L}_{(I\mrm{grad}^{\omega}{\varphi})_{\m{V}}} v_0 \right).
\]
The right-hand side can be written as $-\m{A}^*\m{A}(\varphi)$ following the description \eqref{Eq:linearised_relative_map_nu}. The minus sign follows from the relation
\[
(I\mrm{grad}^{\omega}{\varphi})_{\m{V}} = -\nabla_{\m{V}}\varphi,
\]
where $\nabla_{\m{V}}$ is the vertical Riemannian gradient.
\end{proof}

We define the operator
\begin{equation}\label{Eq:OSC_operator}
\begin{aligned}
\m{G} : \m{P}^{2,l}&\longrightarrow W^{2, l-2}(X) \\
(\varphi, x,v) &\longmapsto p_{E(x,\varphi)} \left(\Theta(\omega, F_\pi (\varphi,\Phi(x)))\right) + \frac{\lambda}{2}\nu_{\varphi, x} \left(F_\pi'(\varphi, \Phi(x),v)\right),
\end{aligned}
\end{equation}
where $\m{U}^{2,l}$ is the space defined in Definition \ref{Def:HandG}, but the functions are considered to be in the Sobolev space $W^{2,l}(X)$ instead of smooth.

\begin{proposition}\label{Prop:openness_OSC}
Let $\pi_X: X \to B$ be a holomorphic submersion with a fibrewise cscK structure $(\omega, I)$ and let $\pi_Y:Y\to B$ be a deformation of $X\to B$ with complex structure $J$. Let $V_\pi$ be the Kuranishi space based at $I$ and let $v_0 \in V_\pi$ represent the complex structure $J$. Assume that
\begin{enumerate}[label = (\roman*)]
\item $(\omega, J)$ is an optimal symplectic connection;
\item the relative automorphism group
\[
\Aut(\pi_Y) := \set*{f \in \Aut(Y, H_Y) \mid f \circ \pi = \pi}
\]
is discrete.
\item the group $\Aut(X_b, H_b)$ is independent of $B$, and will be denoted $G$.
\end{enumerate}
Then for any small deformation $v$ of $v_0$ in $V_\pi^+$ there exists a pair $(x, v) \in TV_\pi$ such that $(\omega, \Phi(x))$ is relatively cscK and $v$ generates a complex structure $J'$, and a K\"ahler potential $\varphi$ such that 
\[
\omega+ i \del\delbar\varphi
\]
is an optimal symplectic connection with respect to $J'$, where the $\del$, $\delbar$ operators are with respect to $\Phi(x)$.
\end{proposition}
\begin{proof}
The proof consists of proving that the operator \ref{Eq:OSC_operator} is an elliptic operator with a trivial kernel so that we can apply the implicit function theorem.
We note that this can be done even though $V_\pi^+$ may be a singular complex space. Indeed, let $\m{X} \to B\times V_\pi^+$ be a family of holomorphic submersions such that the fibre over $0 \in V_\pi^+$ is $X \to B$. By the Kuranishi theorem \cite[\S1]{Kuranishi_family_cpx_str} we can locally consider a smooth trivialisation of the family over $V_\pi^+$ such that the complex structures of the fibrations $\m{X}_x \to B \times \{x\}$ form a smooth family $\{J_s\}$.

As before, let $I$ denote the complex structure of $X$. As we assume that $(\omega, I, v_0)$ is an optimal symplectic connection, $\m{G}(0,0, v_0) = 0$. The derivative with respect to the first component, given by Lemma \ref{Lemma:Derivative_mG1} and \ref{Lemma:Derivative_mG2}, is
\begin{equation*}
\dd_1 \m{G}_{(0,0, v_0)} (\varphi) = -\m{R}^*\m{R}(\varphi) - \lambda\m{A}^*\m{A}(\varphi) = -\widehat{\m{L}}(\varphi).
\end{equation*}
The hypothesis on the automorphism group implies that the kernel of the linearisation is empty, so the implicit function theorem guarantees that there exists a path
\[
(\varphi(x,v), x, v) \in \m{U}^{2,l}
\]
such that locally around $(0,0,v_0)$
\begin{equation}\label{Eq:OSC_operation_form}
\m{G}(\varphi(x,v), x, v) = 0.
\end{equation}
The function $\varphi(x,v)$ is smooth by the standard theory of regularity of solutions to elliptic partial differential equations, applied to the operator $\m{G}$ \cite[Theorem 41]{Besse_EinsteinManifolds}.
Therefore solutions to \eqref{Eq:OSC_operation_form} produce optimal symplectic connections.
\end{proof}

\subsection{The moduli space of optimal symplectic connections}
Let $(Y, H_Y) \to B$ be a relatively K-semistable holomorphic submersion with a degeneration to a relatively cscK fibration $(X, H_X)\to B$. Assume that we have a relatively K\"ahler metric $(\omega, J)$ on $Y$ that is an optimal symplectic connection.
Let $W$ be the subset of the Kuranishi space $V_\pi$ that corresponds to fibrations satisfying the hypotheses of Proposition \ref{Prop:openness_OSC}. Then $W$ is an open subset of the locally closed subvariety $V_\pi^+$ described in \S \ref{Subsec:openness_setting} by Proposition \ref{Prop:openness_OSC}.
\begin{lemma}[{\cite[Corollary to Proposition 2]{Fujiki_CoarseModuliSpace}}]
The group $\Aut(\pi_Y)$ is a subgroup of $\Aut(Y, H_Y)$ with finitely many connected components.
\end{lemma}
In particular, under our assumption $\Aut(\pi_Y)$ is a finite discrete group.
Let $V_{\pi_Y}$ be the Kuranishi space of the fibration $\pi_Y$ and let
\[
\tau : V_{\pi_Y} \to V_\pi
\]
be the map given by completeness of the Kuranishi space. If we denote $\tau^{-1}W =: W_Y$, then $W_Y$ is a locally closed subvariety of $V_{\pi_Y}$.

Let $\m{Y} \to B \times W_Y$ be the Kuranishi family of fibrations which admit an optimal symplectic connection, with central fibration $Y \to B$.
The quotient
\begin{equation}\label{Eq:local_charts}
W_Y / \Aut(\pi_Y)
\end{equation}
is a local complex space and it is Hausdorff since it is the quotient of a variety by a finite group.
We now explain that we can glue the quotients \eqref{Eq:local_charts} to obtain a global Hausdorff moduli space $\m{M}$ of fibrations that admit an optimal symplectic connection.
\begin{remark}
The moduli space $\m{M}$ depends on the group $ G= \Aut_0(X_b, H_b)$. In other words, $\m{M}$ parametrises all fibrations $\pi_Y:Y \to B$ such that they have discrete relative automorphism group and such that they degenerate to a relatively cscK fibration whose fibres have $G$ as their automorphism group.
\end{remark}

The following result builds on \cite[Proposition 6.5]{FujikiSchumacher_Moduli_cscK} and \cite[Lemma 3.8]{FujikiSchumacher_ModuliSymplectic}.

\begin{lemma}\label{Lemma:properness}
Let $\m{Y}$ and $\m{Y}'$ over $W_Y$ be two families of fibrations that admit an optimal symplectic connection.
The group of isomorphisms between $\m{Y}$ and $\m{Y}'$ that preserve the fibration structure, denoted $\mrm{Isom}_{W_Y}(\m{Y}, \m{Y}', B)$, is proper over $W_Y$.
\end{lemma}
\begin{proof}
Let $y_t \to \bar{y}$ be a convergent sequence in $W_Y$ and consider a family of isomorphisms $f_t: \mathcal{Y}_{y_t} \to \mathcal{Y}'_{y_t}$
preserving the projection onto $B$.
Such isomorphisms are fibrewise isometries with respect to the underlying fibrewise Riemannian metrics. Therefore there exists a subsequence $\{f_{t_k}\}$ which converges in the $C^m$-topology to a fibrewise $C^m$-isometry $f: \mathcal{Y}_{\bar{y}} \to \mathcal{Y}'_{\bar{y}}$ \cite[Lemma 3.8]{FujikiSchumacher_ModuliSymplectic}. Moreover, $f$ is a biholomorphism because it is the limit of biholomorphic maps.
Therefore the convergence takes place in $\mrm{Isom}_{W_Y}(\m{Y}, \m{Y}', B)$, which is then proper over $W_Y$.
\end{proof}

The global definition of $\m{M}$ relies on the following lemma, whose proof follows from \cite[XI.6]{ArbarelloCornalbaGriffiths}.

\begin{lemma}\label{Lemma:group_equivalence}
Let $\m{Y} \to B \times W_Y$ be a family of fibrations that admit an optimal symplectic connection. For any points $y,w \in W_Y$, the fibrations $Y_y \to B$ and $Y_w \to B$ are isomorphic if and only if there exists an element $g \in \Aut(\pi_Y)$ such that $g(w) = y$.
\end{lemma}
\begin{proof}
We can prove the lemma for the entire relative Kuranishi space $V_{\pi_Y}$.
The statement follows from:
\begin{enumerate}[label=(\roman*)]
\item For any $y\in V_{\pi_Y}$, the automorphism group $\Aut(\pi_{Y_y})$ is contained in $\Aut(\pi_Y)$;
\item For any $y\in V_{\pi_Y}$ there exists an $\Aut(\pi_{Y_y})$-invariant open neighbourhood $\m{U}_y$ such that any isomorphism between fibres of $\m{Y}\vert_{\m{U}_y}\to\m{U}_y$ is induced by an element of $\Aut(\pi_{Y_y})$.
\end{enumerate}
We begin by proving (ii).
Assume by contradiction that there exist two sequences $\{y_n\}$ and $\{w_n\}$ both converging to $y$. Assume that there exist $g_n$ in $\Aut(\pi_{Y_y}) \setminus\Aut(\pi_Y)$ such that $g_n\cdot w_n =y_n$.
Then by Lemma \ref{Lemma:properness} there exists $g \in \Aut(\pi_{Y_y})$ such that $g_n \to g$.
Up to replacing $g_n$ with $g_n g^{-1}$ and $w_n$ with $g^{-1}w_n$ we may assume that $g$ is the identity of $\Aut(\pi_{Y_y})$.
Consider the map
\[
\begin{aligned}
\Aut(\pi_Y)\times V_{\pi_Y} &\to \widetilde{\mrm{H}}^1_{\m{V}}\\
(\eta,y) &\mapsto \eta\cdot y.
\end{aligned}
\]
By Theorem \ref{Thm:relative_Kuranishi}, this is a local biholomorphism at $(\mrm{id},y)$.
So we have $F(\mrm{id},y_n) = F(g_n,w_n)$. Hence $g_n = \mrm{id}$, a contradiction.

The claim (i) then follows exactly as in \cite[p.204]{ArbarelloCornalbaGriffiths}. We report the proof for completeness. Set
\[
I = \set*{(y,g) \in W_Y \times \Aut(\pi_{y_y})}.
\]
Then $I$ is equal to $\mrm{Isom}(\m{Y},\m{Y},B)$.
Up to shrinking $W_Y$ we can assume that any connected component $I'$ of $I$ intersects $\{0\}\times \mrm{Aut}(\pi_Y)$.
Let $\widehat{I} = \set*{(y,g) \in I' \mid g \in \mrm{Aut}(\pi_Y)}$.
Then $\widehat{I}$ is non-empty and Zariski-closed. It follows from (ii) that $\widehat{I}$ contains an open set. Then $\widehat{I}=I'$, which concludes the proof.
\end{proof}

To glue the charts \eqref{Eq:local_charts}, we use the completeness of the Kuranishi space.
Let $Y_1 \to B$ be another relatively K-semistable fibration which admits an optimal symplectic connection and is close to $Y\to B$. Then the Kurnanishi theorem gives a map
\[
\tau: V_{\pi_{Y_1}} \to V_{\pi_Y}
\]
such that a family $\m{Y}_1 \to B \times V_{\pi_{Y_1}}$ is isomorphic to the pull-back via $\tau$ of $\m{Y} \to B\times V_{\pi_Y}$.
Consider the composition of this map with the inclusion $i : W_{Y_1} \to V_{\pi_{Y_1}}$:
\[
\tau \circ i : W_{Y_1} \to V_{\pi_Y}.
\]
Let $y_1 \in W_{Y_1}$. Then $y_1$ represents a fibration with an optimal symplectic connection, i.e. the image of $y_1$ via the Kuranishi map is a complex structure $J_{y_1}$ such that $(\omega, J_{y_1})$ is an optimal symplectic connection. Therefore, there is a representative $y$ of  $J_{y_1}$ in $V_{\pi_Y}$ which belongs to $W_Y$, so
\begin{equation*}
\alpha:=\tau\circ i : W_{Y_1} \to W_{Y}.
\end{equation*}
Lemma \ref{Lemma:group_equivalence} allows us to pass to the quotient and obtain an isomorphism
\begin{equation}\label{Eq:completeness_nearby_fibres_map_W}
\widetilde{\alpha} : W_{Y_1}/\Aut({\pi_{Y_1}}) \to W_Y/\Aut(\pi_Y),
\end{equation}
which is uniquely determined (while $\alpha$ itself might not be). Indeed, the inverse is constructed by reverting the roles of $Y$ and $Y_1$.
Therefore, we can use it to glue the local charts to give $\m{M}$ the structure of a complex space.


\begin{proposition}
The space $\m{M}$ is a Hausdorff complex space with at most countably many connected components.
\end{proposition}
\begin{proof}
The countability follows from \cite[Theorem 7.3]{FujikiSchumacher_ModuliSymplectic}. We prove the Hausdorff property.
Let $\m{Y} \to B \times W_Y$ and $\m{Y}_1\to B\times W_{Y_1}$ be two families of fibrations which admit an optimal symplectic connection.
Let $y_t \to \bar{y}$ be a sequence in $W_Y$ and $y_{1t} \to \bar{y}_1$ be a sequence in $W_{Y_1}$ and assume that $Y_{y_t}$ is isomorphic to $Y_{1,y_{1t}}$ as fibrations over $B$.
Following the proof of \cite[Proposition 10]{Fujiki_CoarseModuliSpace}, we show that $Y_{\bar{y}}\to B$ is isomorphic to $Y_{1,\bar{y}_1}$.
Let 
\[
\widehat{W} = W_Y \times W_{Y_1},
\]
and let $\widehat{\m{Y}} \to B\times \widehat{W}$ and $\widehat{\m{Y}}_1 \to B \times \widehat{W}$ be the pull-back of $\m{Y}$ and $\m{Y}_1$ using the projections of $\widehat{W}$ onto the first and second component respectively.
Consider
\[
\Sigma = \set{ (y, y_1) \in \widehat{W} \mid Y_y \to B \ \text{is isomorphic to} \ Y_{1, y_1} \to B}.
\]
It follows from the properness of $\mrm{Isom}(\widehat{\m{Y}}, \widehat{\m{Y}}_1, B)$ that $\Sigma$ is a locally closed analytic subvariety of $\widehat{W}$.
Therefore $(\bar{y}, \bar{y}_1) \in \Sigma$, which concludes the proof.
\end{proof}

We have proven the following.
\begin{corollary}
There exists a Hausdorff complex space $\m{M}$ which parametrises holomorphic submersions over a fixed base admitting an optimal symplectic connection, with fixed relative automorphism group. 
\end{corollary}

\section{A Weil-Petersson type K\"ahler metric}\label{Sec:WPmetric}
In this section, we define a K\"ahler metric on the moduli space of fibrations admitting an optimal symplectic connection. We do so by describing a relative version of the theory of Weil-Petersson type metrics developed by Fujiki and Schumacher for cscK manifolds \cite[Sections 8, 9]{FujikiSchumacher_Moduli_cscK}.
In particular, we first define the Weil-Petersson metric locally on the
 sets $W_Y$ and subsequently we extend the definition to the charts $W_Y/\Aut(\pi_Y)$ and then to the moduli space $\m{M}$.
 
We start by recalling that the notions of a smooth Hermitian metric and of a smooth K\"ahler metric are well posed on a complex space {\cite[Definitions 1.1, 1.2]{FujikiSchumacher_Moduli_cscK}}.
Consider a family of holomorphic submersions that admit an optimal symplectic connection, denoted by $\m{Y} \to B \times W_Y$, and let $Y\to B$ be the central fibration.
Let $\omega_t$ be also the optimal symplectic connection on $Y_t \to B$.
For each $k$ sufficiently large, $\omega_{t,k} = \omega_t + k\omega_B$ is a K\"ahler form on $Y_t$, and we denote by $\omega_{t,F}$ its purely vertical part.
Let $\m{Y} \to W_Y$ be the composition with the second projection.

From Theorem \ref{Thm:relative_Kuranishi}, for each $t \in W_Y$ there is an injective map
\begin{equation}\label{Eq:Kuranishi_tangent}
\dd_t\Phi:T_tW_Y \hookrightarrow \widetilde{H}^1_{\m{V}}(Y_t) \subseteq \Omega^{0,1}(\m{V}^{1,0}_{Y_t})
\end{equation}
that identifies a vector $\alpha$ in the Zariski tangent space $T_tW_Y$ with a $(0,1)$-form valued in the $(1,0)$-tangent bundle of $Y_t$, which we will also denote by $\alpha$. The map \eqref{Eq:Kuranishi_tangent} is equivariant with respect to the action of $\Aut(\pi_Y)$.
Therefore we can define an inner product on $T_tW_Y$ by pulling back the $L^2$-product on $\Omega^{0,1}(\m{V}^{1,0}_{Y_t})$ induced by the Hermitian metric associated to $\omega_{t,F} + \omega_B$. For any $\alpha, \beta \in T_tW_Y$, its imaginary part is given by
\begin{equation}\label{Eq:WP_metric}
\Omega_t(\alpha,\beta) := \langle\alpha,\beta\rangle_{\omega_{t,F}+\omega_B} = \int_{Y_t} \Lambda_{\omega_{t,F}+\omega_B}\mrm{Tr}_{\omega_{t,F}} (\alpha\overline{\beta}) \omega_t^m\wedge\omega_B^n,
\end{equation}
where we denote by $\Lambda$ the contraction of the covariant part and by $\mrm{Tr}$ the trace of the contravariant part.
We give the following definition.
\begin{definition}\label{Def:WP_metric}
The \emph{relative Weil-Petersson metric} on $W_Y$, denoted by $\Omega_{WP}$, is the two-form $\{\Omega_{t}\}_{t \in W_Y}$.
\end{definition}
Using the compatibility of the deformations with the K\"ahler form, we write the trace in coordinates as
\[
\begin{aligned}
\Lambda_{\omega_{t,F}+\omega_B}\mrm{Tr}_{\omega_{t,F}} (\alpha\overline{\beta}) = \alpha^a_{\ \bar{q}} \overline{\beta^b_{\ \bar{p}}} (\omega_{t,F})_{a\bar{b}} (\omega_{t,F}+\omega_B)^{\bar{q}p}= \alpha^a_{\ \bar{b}} \overline{\beta^b_{\ \bar{a}}} + \Lambda_{\omega_B} \mrm{Tr}_{\omega_{t,F}} (\alpha\overline{\beta}).
\end{aligned}
\]
Therefore the integral \eqref{Eq:WP_metric} can be written as the sum
\[
\int_{Y_t} \alpha^a_{\ \bar{b}} \overline{\beta^b_{\ \bar{a}}} \omega_t^m\wedge\omega_B^n + \int_{Y_t} \Lambda_{\omega_B}\mrm{Tr}_{\omega_{t,F}} (\alpha\overline{\beta}) \omega_t^m\wedge\omega_B^n.
\]

\begin{remark}
The first term can be split over $Y_t$ as
\begin{equation}\label{Eq:spit_LOT_integral}
\int_B \left( \int_{Y_{t,b}}\alpha^a_{\ \bar{b}} \overline{\beta^b_{\ \bar{a}}} \omega_t^m \right)\omega_B^n.
\end{equation}
In particular, it vanishes when $\alpha$ and $\beta$ restrict to the trivial deformation on the fibres.
This is the case when $\m{Y}\to B \times W_Y$ is a family of holomorphic submersions with rigid fibres, for example.
\end{remark}

We next prove that \eqref{Eq:WP_metric} is closed and positive definite. Consider for each $t\in W_Y$ a $k$-dependent inner product on $\widetilde{H}^1_{\m{V}}(Y_t)$, defined using the K\"ahler form $\omega_{t,k}$ of $Y_t$:
\[
\Omega_{k,t}(\alpha,\beta) := \int_{Y_t} \langle \alpha, \beta\rangle_{\omega_{t,k}} \omega_{t,k}^{n+m} = \int_{Y_t} \Lambda_{\omega_{t,k}}\mrm{Tr}_{\omega_{t,k}} (\alpha\overline{\beta}) \omega_{t,k}^{n+m}.
\]
The collection $\{\Omega_{t,k}\} =: \Omega_k$ is the Weil-Petersson type Hermitian metric defined by Fujiki-Schumacher \cite[\S7]{FujikiSchumacher_Moduli_cscK} for any family of smooth polarised varieties.
Using Lemma \ref{Lemma:contraction_1-1_tensor}, we can write in local holomorphic coordinates
\[
\begin{aligned}
\langle \alpha, \beta\rangle_{\omega_{t,k}} = \alpha^a_{\ \bar{q}} \overline{\beta^b_{\ \bar{p}}} (\omega_{t,k})_{a\bar{b}} (\omega_{t,k})^{\bar{q}p}= \alpha^a_{\ \bar{b}} \overline{\beta^b_{\ \bar{a}}} + k^{-1}\Lambda_{\omega_B} \mrm{Tr}_{\omega_{t,F}} (\alpha\overline{\beta}) + \Ok{-2}.
\end{aligned}
\]
The expansion in powers of $k$ of $\omega_{t,k}$ reads
\[
\omega_{t,k}^{n+m} = k^n\omega_t^m \wedge \omega_B^n + \frac{k^{n-1}}{nm}\omega_t^{m+1}\wedge \omega_B^{n-1} + \Ok{n-2}.
\]
Then
\begin{equation}\label{Eq:WP_k}
\begin{split}
\Omega_{t,k}(\alpha,\beta)= k^n\int_{Y_t}\alpha^a_{\ \bar{b}} \overline{\beta^b_{\ \bar{a}}} \omega_t^m\wedge\omega_B^n+k^{n-1}\left[\int_{Y_t} \alpha^a_{\ \bar{b}} \overline{\beta^b_{\ \bar{a}}}\omega_t^{m+1}\wedge\omega_B^{n-1}\right.\\
\left. + \int_{Y_t}\Lambda_{\omega_B} \mrm{Tr}_{\omega_{t,F}} (\alpha\overline{\beta}) \omega_t^m \wedge \omega_B^n\right] +\Ok{n-2}.
\end{split}
\end{equation}
Then the two terms in the sum \eqref{Eq:WP_metric} are the first and third coefficients of this expansion.

We next describe a fibre integral formula for the Weil-Petersson metric on $W_Y$.
Let $\omega_{\m{Y},k}$ be the relatively K\"ahler metric on $W_Y$ such that its restriction to each $Y_t$ is the K\"ahler metric $\omega_t + k\omega_B$, where $\omega_t$ is an optimal symplectic connection on $Y_t \to B$. 
Let also $\rho_{\m{Y},k}$ be the curvature of the Hermitian structure induced by $\omega_{\m{Y},k}$ on the relative anticanonical bundle
\[
-K_{\m{Y}/W_Y} = \ext{n+m} \m{V}_{\m{Y}/W_Y},
\]
where $\m{V}_{\m{Y}/W_Y}$ denotes the vertical tangent bundle of the fibration $\m{Y}\to W_Y$.
We first recall the following result.

\begin{proposition}[{\cite[\S8]{FujikiSchumacher_Moduli_cscK}}]\label{Prop:fibreintegral}
The $k$-dependent $(1,1)$-form $\Omega_k$ can be written as a fibre integral over the map $\m{Y} \to W_Y$:
\begin{equation}\label{Eq:fibre_integral_k}
\Omega_{k}(\omega_{\m{Y},k})=\int_{\m{Y}/W_Y} \rho_{\m{Y},k} \wedge \omega_{\m{Y},k}^{n+m} - \frac{1}{n+m+1}\int_{\m{Y}/W_Y} \mrm{Scal}_{\m{V}}(\omega_{\m{Y},k}) \omega_{\m{Y},k}^{n+m+1}.
\end{equation}
\end{proposition}
By expanding $\Omega_{k}(\omega_{\m{Y},k})$ in powers of $k$, we can find a fibre integral formula for the Weil-Petersson metric \eqref{Eq:WP_metric}.
Since the base $B$ is fixed, the relative metric $\omega_{\m{Y},k}$ can be written as
\[
\omega_{\m{Y},k} = \widehat{\omega} + k\omega_B,
\]
were the restriction of $\widehat{\omega}$ to each $Y_t$ is the optimal symplectic connection $\omega_t$.
Then $\rho_{\m{Y},k}$ can be expanded in powers of $k$ as
\[
\begin{aligned}
\rho_{\m{Y},k} &= i\del\delbar\log \left(k^n\widehat{\omega}^m\wedge\omega_B^n + \Ok{n-1}\right)\\
&= i\del\delbar\log \left(\widehat{\omega}^m\wedge\omega_B^n\right) + \Ok{-1}\\
&= i\del\delbar \log \det(\widehat{\omega}) + i \del \delbar \log\det(\omega_B) + \Ok{-1},
\end{aligned}
\]
where the second line follows from the fact that $\widehat{\omega}^m\wedge\omega_B^n$ is a volume form. The $k^{-1}$-term is exact, because the two volume forms $\omega_{\m{Y},k}^{m+n}$ and $\widehat{\omega}^m\wedge\omega_B^n$ both induce the class $c_1(-K_{\m{Y}/W_Y})$.
Moreover, $\det(\widehat{\omega})$ is the relative determinant of $\widehat{\omega}$ and $i\del\delbar \log \det(\widehat{\omega})$ is the curvature of the Hermitian metric induced by $\widehat{\omega}$ on the relative anticanonical bundle $-K_{\m{Y}/B\times W_Y}$.
To compute the expansion of the vertical scalar curvature, we use the expression \cite[Proposition 5.4]{Ortu_OSCdeformations} for the scalar curvature of $\omega_k$, which we report here.

\begin{proposition}\label{Prop:expansion_scalar_curvature}
The scalar curvature of $\omega_k$ on the submersion $Y\to B$ admits an expansion
\begin{equation*}
\mrm{Scal}(\omega_k) = \widehat{S}_b + k^{-1}\left( \mrm{Scal}(\omega_B) - \Lambda_{\omega_B}\alpha + p_E( \Theta(\omega))\right) + O\left(k^{-3/2}\right)
\end{equation*}
where:
\begin{enumerate}[label=(\roman*)]
\item $\widehat{S}_b$ is the average scalar curvature of the fibres;
\item $\alpha$ is the Weil-Petersson metric on $B$ induced by the relatively cscK degeneration $X \to B$.
\end{enumerate}
\end{proposition}

Since $\widehat{\omega}$ is an optimal symplectic connection when restricted to each $Y_t$, using the proposition \ref{Prop:expansion_scalar_curvature} the vertical scalar curvature of $\omega_{\m{Y},k}$ admits an expansion as
\[
\mrm{Scal}_{\m{V}}(\omega_{\m{Y},k}) = \widehat{S}_b + k^{-1}\left( \mrm{Scal}(\omega_B)- \Lambda_{\omega_B}\widehat{\alpha}\right) + \Ok{-2},
\]
where $\widehat{\alpha}$ is a closed two-form on $B$.
Then the leading order term of \eqref{Eq:fibre_integral_k} is
\[
I_0=\int_{\m{Y}/W_Y} i\del\delbar \log \det(\widehat{\omega})\wedge\widehat{\omega}^{m}\wedge\omega_B^{n} - \frac{1}{n+m+1}\int_{\m{Y}/W_Y} \widehat{S}_b \widehat{\omega}^{m+1}\wedge\omega_B^{n}
\]
The sub-leading order term is the sum of the four integrals
\begin{align}
&I_1=\int_{\m{Y}/W_Y} i\del\delbar \log \det(\widehat{\omega})\wedge\widehat{\omega}^{m+1}\wedge\omega_B^{n-1}, \label{Eq:I_1} \\
&I_2 = - \frac{1}{n+m+1}\int_{\m{Y}/W_Y} \widehat{S}_b \widehat{\omega}^{m+2}\wedge\omega_B^{n}, \notag \\
&I_3 =  \frac{1}{n}\int_{\m{Y}/W_Y} \mrm{Scal}(\omega_B)\widehat{\omega}^{m+1}\wedge\omega_B^n, \notag \\
&I_4 = - \frac{1}{n+m+1}\int_{\m{Y}/W_Y} \left(\mrm{Scal}(\omega_B)- \Lambda_{\omega_B}\widehat{\alpha}\right)  \widehat{\omega}^{m+1}\wedge\omega_B^n. \notag
\end{align}

We can use this expansion to prove a fibre integral formula in our setting.
\begin{lemma}\label{Lemma:fibre_integral}
The Weil-Petersson metric $\Omega_{WP}(\omega_{\m{Y}})$ \eqref{Eq:WP_metric} can be written as the fibre integral
\begin{equation}\label{Eq:fibre_integral}
I_0 +I_2 + I_3 + I_4.
\end{equation}
\end{lemma}
\begin{proof}
It follows from \cite[Lemma 8.5]{FujikiSchumacher_Moduli_cscK} applied to the family $\m{Y}\to B\times W_Y$ that the two-form defined as the collection of the integrals
\[
\int_{Y_t} \alpha^a_{\ \bar{b}} \overline{\beta^b_{\ \bar{a}}}\omega_t^{m+1}\wedge\omega_B^{n-1}
\]
is equal to $I_1$. Indeed, working locally in $D \subset B\times W_Y$, Fujiki and Schumacher prove that one can trivialise the family $\m{Y}$ over $D$ as $Y_{0,0}\times D$ such that the horizontal distribution induced by $\widehat{\omega}$ is preserved.
Then, given a family $\{\beta_t\}$ of vertical deformations of the complex structure of $Y_{0,0}$ which represent the family $Y_{0,0}\times D \to D$, \cite[Lemma 8.5]{FujikiSchumacher_Moduli_cscK} gives the equality
\[
i\del\delbar \log \det(\widehat{\omega}) = \mrm{Tr} (\del_t\beta_t\vert_{t=0}\overline{\del_t\beta_t}\vert_{t=0}),
\]
where $\del_t\beta_t$ is the map \eqref{Eq:Kuranishi_tangent}.
\end{proof}

\begin{lemma}\label{Lemma:WP_closed_positive}
The two-form $\Omega_{WP}$ is closed and positive-definite on $W_Y$.
\end{lemma}
\begin{proof}
Since the map \eqref{Eq:Kuranishi_tangent} is injective, the integral \eqref{Eq:WP_metric} is positive.
To prove closedness, we show that the terms $I_0$, $I_2$, $I_3$ and $I_4$ in Lemma \ref{Lemma:fibre_integral} are closed.
The terms $I_0$ and $I_2$ are closed because they are the fibre integrals of a closed form.
The term $I_4$ can be written as
\[
\int_{\m{Y}/W_Y} (\rho_B - \widehat{\alpha}) \wedge \widehat{\omega}^{m+1} \wedge \omega_B^{n-1},
\]
where $\rho_B$ is the Ricci form of $\omega_B$. In particular $(\rho_B - \widehat{\alpha})\wedge\omega_B^{n-1}$ is a top degree form on $B$, hence it is closed. So $I_4$ is closed. Analogously, $I_3$ is closed.
\end{proof}

Let $h_{WP}(W_Y)$ be the Hermitian metric on the tangent bundle to $W_Y$ induced by the two-form $\Omega_{WP}$.
\begin{theorem}
The Hermitian metric $h_{WP}(W_Y)$ induces a global K\"ahler metric on the moduli space $\m{M}$.
\end{theorem}
\begin{proof}
The theorem follows from the fact that the action of the finite group $\Aut(\pi_Y)$ is induced by an automorphism of the Kuranishi family, so the Weil-Petersson metric $h_{WP}(W_Y)$ is invariant for the action of $\Aut(\pi_Y)$.
Therefore it defines a metric on the quotient $W_Y/\Aut(\pi_Y)$.
\end{proof}

\begin{remark}
We have defined a Weil-Petersson metric on $W_Y$ that is independent of the adiabatic parameter $k$.
This is reasonable from the point of view of describing the moduli space of fibrations using the optimal symplectic connection alone, which is only relatively K\"ahler.
However, it is possible that a different kind of a Weil-Petersson type metric could be defined by taking a sequence of K\"ahler metrics that depend on $k$, where the adiabatic construction plays a bigger role.
\end{remark}

\subsection{The determinant line bundle for the Weil-Petersson metric}
We construct a line bundle on $W_Y$, and hence on the moduli space, such that the Weil-Petersson metric represents its first Chern class.
To do so, we appeal to the theory of Deligne pairings \cite{Elkik_fibres_integral_Chernclasses}, \cite{Elkik_metrics_Deligne_pairings}, \cite[\S1]{BoucksomHisamotoJonsson_DelignePairings}.
Let $M \to B$ be a flat, projective morphism between complex algebraic varieties of relative dimension $d$ and consider $d+1$ line bundles $L_0,..., L_d$ on $M$. The push-forward of the intersection product of $L_0, ..., L_d$ is an isomorphism class of line bundles on $B$, represented by the cohomology class
\begin{equation}\label{Eq:Deligne_pairing_class}
\int_{M/B} c_1(L_0) \wedge \dots \wedge c_1(L_d).
\end{equation}
The \emph{Deligne pairing} of $L_0, \dots, L_d$, denoted by $
\langle L_0 , \dots, L_d \rangle_{M/B} $, is a canonical choice of a line bundle on $B$ such that \eqref{Eq:Deligne_pairing_class} is its first Chern class.
The construction is symmetric, multilinear and functorial.
Moreover, if $h_0, \dots, h_d$ are Hermitian metrics on $L_0, \dots, L_d$ respectively, the theory provides a metric $\langle h_0, \dots, h_d\rangle_{M/B}$ on $\langle L_0 , \dots, L_d \rangle_{M/B} $.
Denoting by $\omega_0, \dots, \omega_d$ the curvature forms of $h_0, \dots, h_d$ respectively, the curvature of $\langle h_0, \dots, h_d\rangle_{M/B}$ is given by the fibre integral
\begin{equation}\label{Eq:fibre_integral_general}
\int_{M/B} \omega_0 \wedge \dots \wedge \omega_d.
\end{equation}

The fibre integral formula for the Weil-Petersson metric on $W_Y$ of Lemma \ref{Lemma:fibre_integral} is a special case of the expression \eqref{Eq:fibre_integral_general}.
To describe it as the curvature form of a line bundle on $W_Y$,
we first recall the following result of Fujiki and Schumacher for the $k$-dependent Weil-Petersson metric.

\begin{proposition}[{\cite[\S9]{FujikiSchumacher_Moduli_cscK}}]\label{Prop:fibreintegral_CMlinebundle_k}
The $k$-dependent Weil-Petersson type K\"ahler metric $\Omega_{k}(\omega_{\m{Y},k})$ represents the first Chern class of the line bundle
\begin{equation}\label{Eq:CM_line_bundle_k}
\langle-K_{\m{Y}/W_Y}, \m{L}_{k}^{n+m}\rangle_{\m{Y}/W_Y} - \frac{1}{n+m+1}\frac{-K_Y\cdot (H_Y+kL)^{n+m-1}}{(H_Y+kL)^{n+m}}\langle \m{L}_{k}^{n+m+1}\rangle_{\m{Y}/W_Y},
\end{equation}
where the constant
\[
\widehat{S}_{Y} := \frac{-K_Y\cdot (H_Y+kL)^{n+m-1}}{(H_Y+kL)^{n+m}}
\]
is the average scalar curvature of $Y$, and hence of each $Y_t$, with respect to the metric $\omega + k\omega_B$.
\end{proposition}

We use Proposition \ref{Prop:fibreintegral_CMlinebundle_k} to  prove the following.
\begin{proposition}
There exists a line bundle $\m{D}(Y)$ on $W_Y$ whose first Chern class is represented by the Weil-Petersson metric \eqref{Eq:WP_metric}.
\end{proposition}
\begin{proof}

Let
$
\widehat{\m{H}} \to \m{Y}
$
be the relatively ample line bundle induced by each relative polarisation $H_t \to Y_t$. More precisely, we consider the fibration $\m{Y} \to W_Y$ as the composition of
\[
\m{Y} \to B \times W_Y \to W_Y
\]
so $\widehat{\m{H}}$ is relatively ample with respect to the first projection.
Moreover, since the base $B$ is fixed, we can pull-back $L$ to $\m{Y}$ and consider it as a line bundle on $\m{Y}$.
So we can define a relatively ample line bundle $\m{L}_k$ over $\m{Y}$ as
\begin{equation}\label{Eq:relative_polarisation}
\m{L}_k = \widehat{\m{H}} \otimes kL,
\end{equation}
whose fibre over $Y_t$ is $\m{L}_k\vert_{t} = H_t \otimes kL$.
Its first Chern class contains the relative metric $\omega_{\m{Y},k}$, where $\widehat{\omega}$ is in $c_1(\widehat{\m{H}})$.
We define a line bundle $\m{D}(Y)\to W_Y$ by using the fibre integral formula \eqref{Eq:fibre_integral} and the expansion in powers of $k$ of the line bundle \eqref{Eq:CM_line_bundle_k}.
Expanding in $k$ the intersection product $\m{L}_k^{n+m}$, $\m{L}_k^{n+m+1}$ and the expression \eqref{Eq:CM_line_bundle_k} we obtain that the Weil-Petersson metric $\Omega_{WP}(\omega_{\m{Y},k})$ represents the first Chern class of the line bundle
$\m{D}(Y)$ given as the tensor product of the line bundles
\[
\begin{aligned}
&\m{D}_0(Y) = \langle-K_{\m{Y}/B\times W_Y}, \widehat{\m{H}}^{m}, L^{n} \rangle_{\m{Y}/{W_Y}}, \\
&\m{D}_2(Y)=- \frac{1}{n+m+1}\frac{-K_{Y/B}\cdot L^{n}\cdot H_Y^{m-1}}{L^n\cdot H_Y^m}\left(\langle  \widehat{\m{H}}^{m+1}, L^{n} \rangle_{\m{Y}/{W_Y}} + \langle  \widehat{\m{H}}^{m+2}, L^{n-1} \rangle_{\m{Y}/{W_Y}}\right), \\
&\m{D}_3(Y)=\langle -K_B,  \widehat{\m{H}}^{m+1}, L^{n-1} \rangle_{\m{Y}/{W_Y}}, \\
&\m{D}_4(Y)= - \frac{1}{n+m+1} \frac{\left(-K_{Y/B} \cdot H_Y^{m-1}\cdot L^n\right)\left(H_Y^{m+1}\cdot L^{n-1}\right)}{H_Y^m\cdot L^n}\langle  \widehat{\m{H}}^{m+1}, L^{n} \rangle_{\m{Y}/{W_Y}},
\end{aligned}
\]
defined using the Deligne pairing.
Indeed, by expanding the expression \eqref{Eq:CM_line_bundle_k} in powers of $k$ we see that the sum of the leading order term and the sub-leading order term is given by
\[
\m{D}_0(Y)+\m{D}_1(Y)+\m{D}_2(Y)+\m{D}_3(Y)+\m{D}_4(Y),
\]
where $\m{D}_1(Y)$ is given by
\[
\m{D}_1(Y) = \langle-K_{\m{Y}/B\times W_Y}, \widehat{\m{H}}^{m+1}, L^{n-1} \rangle_{\m{Y}/{W_Y}}.
\]
However, its first Chern class is represented by the term \eqref{Eq:I_1}, which does not appear in the fibre integral formula for the relative Weil-Petersson metric of Lemma \ref{Lemma:fibre_integral}.
This concludes the proof.
\end{proof}

We have constructed a line bundle $\m{D}_Y$ on $W_Y$ whose first Chern class is the Weil-Petersson metric $\Omega_{WP}(\omega_{\m{Y},k})$.
Let now $p : W_Y \to \m{M}$ be the composition of the maps
\[
W_Y \to W_Y/\Aut(\pi_Y) \hookrightarrow \m{M},
\]
where the first map is the quotient by the group action and the second map is the inclusion of a local chart in the moduli space.
It follows from Lemma \ref{Lemma:group_equivalence} that the orbits of the $\Aut(\pi_Y)$-action on $W_Y$ correspond to isomorphic manifolds in the family $\m{Y} \to W_Y$.
Therefore the line bundle $\m{D}(Y)$ is invariant for the action of $\Aut(\pi_Y)$ and thus descends to a line bundle $\widetilde{\m{D}}(Y)$ on the quotient $W_Y/\Aut(\pi_Y)$.

\begin{theorem}\label{Thm:CM_line_bundle}
There exists a line bundle $\m{F}$ on the moduli space $\m{M}$ such that its restriction to each chart $W_Y/\Aut(\pi_Y)$ it is isomorphic to $\widetilde{\m{D}}(Y)$. 
\end{theorem}
\begin{proof}
Let $W_{Y_1}/\Aut(\pi_{Y_1})$ and $W_{Y}/\Aut(\pi_{Y})$ be two local charts of $\m{M}$ with non empty intersection.
Then using completeness of the Kuranishi space there exists an isomorphism $\widetilde{\alpha}: W_{Y_1}/\Aut(\pi_{Y_1}) \to W_Y/\Aut(\pi_Y)$ \eqref{Eq:completeness_nearby_fibres_map_W} which preserves the relative polarisation \eqref{Eq:relative_polarisation} and the submersions onto the base $B$.
By functoriality of the Deligne pairings, the pull-back $\widetilde{\alpha}^*\widetilde{\m{D}}(Y_1)$ is then isomorphic to $\widetilde{\m{D}}(Y)$.
Therefore, on the intersection of $W_Y/\Aut(\pi_{Y})$ and $W_{Y_1}/\Aut(\pi_{Y_1})$ there is a morphism of line bundles $\chi : \widetilde{\m{D}}(Y_1) \overset{\sim}{\to} \widetilde{\m{D}}(Y)$.
Let $\varphi_Y : \widetilde{\m{D}}(Y) \to W_Y/\Aut(\pi_{Y}) \times \bb{C}$ and $\varphi_{Y_1} : \widetilde{\m{D}}(Y_1) \to W_{Y_1}/\Aut(\pi_{Y_1}) \times \bb{C}$ be local trivialisations and, on the intersection,
\[
\psi_{Y_1Y} := \varphi_Y\circ\chi\circ\varphi_{Y_1}^{-1}.
\]
The map $\psi_{Y_1Y}$, viewed as a function on $\bb{C}$ is invertible, with inverse $\psi_{YY_1}$.
Indeed, the map $\chi$ is an isomorphism because $\widetilde{\alpha}$ is.
The same argument proves that the cocycle condition holds.
\end{proof}

The following corollary is a consequence of Theorem \ref{Thm:CM_line_bundle} and of \cite[Proposition 1.7]{FujikiSchumacher_Moduli_cscK}.
\begin{corollary}
Any compact analytic subspace of $\m{M}$ is projective.
\end{corollary}

\bibliographystyle{acm} 
\bibliography{../bibliografia}

\end{document}